\documentclass[12pt]{article}

\usepackage{amsmath,amssymb}

\newcommand{\non}{\nonumber}

\newcommand{\wh}{\widehat}
\newcommand{\ot}{\otimes}
\newcommand{\lan}{\langle\ts}
\newcommand{\ran}{\ts\rangle}
\newcommand{\ve}{\varepsilon}
\newcommand{\ts}{\,}

\newcommand{\U}{ {\rm U}}
\newcommand{\Z}{ {\rm Z}}
\newcommand{\Sr}{ {\rm S}}
\newcommand{\Y}{ {\rm Y}}
\newcommand{\C}{ {\rm C}}
\newcommand{\M}{ {\rm M}}

\newcommand{\tr}{ {\rm tr}\ts}

\newcommand{\R}{ {\rm R}}
\newcommand{\J}{ {\rm J}}
\newcommand{\F}{ {\rm F}}
\newcommand{\Norm}{ {\rm Norm}\ts}
\newcommand{\End}{{\rm{End}\ts}}
\newcommand{\Hom}{{\rm{Hom}}}
\newcommand{\sgn}{ {\rm sgn}\ts}

\newcommand{\CC}{\mathbb{C}}
\newcommand{\ZZ}{\mathbb{Z}}

\newcommand{\Sym}{\mathfrak S}
\newcommand{\h}{\mathfrak h}
\newcommand{\gl}{\mathfrak{gl}}
\newcommand{\oa}{\mathfrak{o}}
\newcommand{\spa}{\mathfrak{sp}}
\newcommand{\g}{\mathfrak{g}}
\newcommand{\agot}{\mathfrak{a}}
\newcommand{\pgot}{\mathfrak{p}}
\newcommand{\kgot}{\mathfrak{k}}
\newcommand{\sll}{\mathfrak{sl}}

\newcommand{\la}{\lambda}
\newcommand{\La}{\Lambda}
\newcommand{\xiL}{\xi^{}_{\La}}
\newcommand{\al}{\alpha}
\newcommand{\be}{\beta}

\newcommand{\Proof}{\noindent{\it Proof.}\ \ } 
\newcommand{\Outline}{\noindent{\it Outline of the proof.}\ \ }

\newtheorem{thm}{Theorem}[section]
\newtheorem{prop}[thm]{Proposition}
\newtheorem{lem}[thm]{Lemma}
\newtheorem{cor}[thm]{Corollary}
\newtheorem{defin}[thm]{Definition}
\newtheorem{example}[thm]{Example}
\newtheorem{remark}[thm]{Remark}

\newcommand{\bth}{\begin{thm}}
\renewcommand{\eth}{\end{thm}}
\newcommand{\bpr}{\begin{prop}}
\newcommand{\epr}{\end{prop}}
\newcommand{\ble}{\begin{lem}}
\newcommand{\ele}{\end{lem}}
\newcommand{\bco}{\begin{cor}}
\newcommand{\eco}{\end{cor}}
\newcommand{\bex}{\begin{example}}
\newcommand{\eex}{\end{example}}
\newcommand{\bre}{\begin{remark}}
\newcommand{\ere}{\end{remark}}
\newcommand{\bde}{\begin{defin}}
\newcommand{\ede}{\end{defin}}
\newcommand{\bal}{\begin{aligned}}
\newcommand{\eal}{\end{aligned}}
\newcommand{\beq}{\begin{equation}}
\newcommand{\ben}{\begin{equation*}}
\newcommand{\bsp}{\begin{split}}
\newcommand{\esp}{\end{split}}

\def\beql#1{\begin{equation}\label{#1}}

\newbox\squ  
\setbox\squ=\hbox{\vrule width.3pt
            \vbox{\hrule height.3pt width.4em\kern1ex\hrule height.3pt}%
            \vrule width.3pt}
\def\endproof{%
  \ifmmode\eqno\copy\squ\smallskip\else{\unskip\nobreak\hfil%
    \penalty50\hskip2em\hbox{}\nobreak\hfil\copy\squ
    \parfillskip=0pt \finalhyphendemerits=0\penalty-100\smallskip}
  \fi}

\begin{document}

\title{\Large\bf Gelfand--Tsetlin bases
for classical Lie algebras\thanks{\noindent
To appear in {\it Handbook of Algebra\/} (M. Hazewinkel, Ed.),
Elsevier.
}
}
\author{{\sc A. I. Molev}\\[15pt]
School of Mathematics and Statistics\\
University of Sydney,
NSW 2006, Australia\\
{\tt
alexm\hspace{0.09em}@\hspace{0.1em}maths.usyd.edu.au
}
}

\date{}
\maketitle


\newpage

\tableofcontents

\newpage

\pagestyle{plain}
\setcounter{page}{1}

\section{Introduction}\label{sec:int}
\setcounter{equation}{0}

The theory of semisimple Lie algebras and their representations is in the heart
of modern mathematics. It has numerous connections
with other areas of mathematics and physics.
The simple Lie algebras over the field of complex numbers
were classified in the works of Cartan and Killing in the 1930's.
There are four infinite series $A_n$, $B_n$, $C_n$, $D_n$
which are called the {\it classical Lie algebras\/}, and
five exceptional Lie algebras $E_6$, $E_7$, $E_8$, $F_4$, $G_2$. 
The structure of these Lie algebras is uniformly described
in terms of certain finite sets of vectors
in a Euclidean space called the {\it root systems\/}.
Due to Weyl's complete reducibility theorem, the theory of
finite-dimensional representations of the semisimple Lie algebras
is largely reduced to the study of irreducible representations.
The irreducibles are parametrized by their {\it highest weights\/}. 
The characters and dimensions are explicitly known by
the {\it Weyl formula\/}. The reader is refered to, e.g., the books of
Bourbaki~\cite{b:ga}, Dixmier~\cite{d:ae},
Humphreys~\cite{h:il} or Goodman and Wallach~\cite{gw:ri}
for a detailed exposition of the theory.

However, the Weyl formula for the dimension does not use any explicit
construction of the representations. Such constructions remained unknown
until 1950 when Gelfand and Tsetlin\footnote{Some authors
and translators write this name in English
as {\it Zetlin, Tzetlin, Cetlin, {\rm or} Tseitlin.}}
published two short papers \cite{gt:fdu}
and \cite{gt:fdo} (in Russian) where they solved the problem for the 
general linear Lie algebras	(type $A_n$) and the orthogonal Lie algebras
(types $B_n$ and $D_n$), respectively. Later, Baird and Biedenharn \cite{bb:rs} (1963)
commented on \cite{gt:fdu} as follows: 

\begin{itemize}
\item[]
``This paper is extremely 
brief (three pages)
and does not appear to have been translated in either the usual Journal
translations or the translations on group-theoretical subjects of the American
Mathematical Society, or even referred to in the review articles on group theory
by Gelfand himself. Moreover, the results are presented without the slightest
hint as to the methods employed and contain not a single reference
or citation of other work. In an effort to understand the meaning
of this very impressive work, we were led to develop the proofs~\dots."
\end{itemize}

\noindent
Baird and Biedenharn employed the calculus of {\it Young patterns\/}
to derive the Gelfand--Tsetlin formulas.\footnote{ 
An indication
of the proof of the formulas of \cite{gt:fdu}
is contained in a footnote in the paper
Gelfand and Graev~\cite{gg:fd} (1965). It says that
for the proof one has ``to verify the commutation relations \dots;
this is done by direct calculation".}
Their interest to the formulas
was also motivated by the connection with the fundamental
{\it Wigner coefficients\/}; see Section~\ref{subsec:ci} below.

A year earlier (1962) Zhelobenko published an independent work~\cite{z:cg}
where he derived the {\it branching rules\/}
for all classical Lie algebras.	In his approach 
the representations are realized in a space of
polynomials satisfying the ``indicator system" of differential equations.
He outlined a method to construct the {\it lowering operators\/}
and to derive the matrix element formulas for the 
case of the general linear Lie algebra $\gl_n$.
An explicit ``infinitesimal" form for the lowering operators as elements of 
the enveloping algebra was found by Nagel and Moshinsky~\cite{nm:ol} (1964)
and independently by Hou Pei-yu~\cite{h:ob} (1966). The latter work
relies on Zhelobenko's results \cite{z:cg} and also contains
a derivation of the Gelfand--Tsetlin formulas alternative to
that of Baird and Biedenharn.
This approach was further
developed in the book by Zhelobenko~\cite{z:cl}
which contains its detailed account.

The work of Nagel and Moshinsky
was extended to the orthogonal Lie algebras $\oa_N$
by Pang and Hecht~\cite{ph:lr} and Wong~\cite{w:ro} who 
produced explicit infinitesimal expressions for the lowering operators and
gave a derivation of the formulas of Gelfand and Tsetlin~\cite{gt:fdo}.

During the half a century passed since the work of Gelfand and Tsetlin,
many different approaches were developed to construct bases
of the representations of the classical Lie algebras.
New interpretations of the lowering operators
and new proofs of the Gelfand--Tsetlin formulas 
were discovered by several authors. In particular,
Gould~\cite{g:cir, g:ia, g:me, g:wc} employed
the {\it characteristic identities\/} 
of Bracken and Green~\cite{bg:vo, gr:ci} to calculate
the Wigner coefficients and matrix elements of generators
of $\gl_n$ and $\oa_N$.
The {\it extremal projector\/}
discovered by Asherova, Smirnov and Tolstoy~\cite{ast:po, ast:po2, ast:dc}
turned out to be
a powerful instrument in the representation 
theory of the simple Lie algebras. It plays an essential role
in the theory of {\it Mickelsson algebras\/} developed by Zhelobenko 
which has a wide spectrum
of applications from the branching rules and reduction problems to
the classification of Harish-Chandra modules; see
Zhelobenko's expository paper~\cite{z:it} and his book~\cite{z:rr}.
Two different {\it quantum minor\/} interpretations of the lowering
and raising operators were given by Nazarov and Tarasov~\cite{nt:yg}
and the author~\cite{m:gt}. These techniques are based on the theory
of quantum algebras called the {\it Yangians\/} and allow an independent derivation
of the matrix element formulas.
We shall discuss the above approaches in more detail in Sections~\ref{subsec:ama},
\ref{subsec:ci} and \ref{subsec:qm} below.

A quite different method to construct modules over
the classical Lie algebras is developed in the papers by
King and El-Sharkaway~\cite{ke:sy}, Berele~\cite{b:cs},
King and Welsh~\cite{kw:co},
Koike and Terada~\cite{kt:yd}, Proctor~\cite{p:yt}, Nazarov~\cite{n:rt}. 
In particular,
bases in the representations of the orthogonal and symplectic
Lie algebras parametrized by $\oa_N$-standard
or $\spa_{2n}$-standard Young tableaux are constructed.
This method provides an algorithm for
calculation of the representation matrices.
It is based on the 
Weyl realization of the representations of the classical
groups in tensor spaces; see Weyl~\cite{w:cg}.
A detailed exposition of the theory of the classical groups
together with many recent developments are presented
in the book by Goodman and Wallach~\cite{gw:ri}.

Bases with special properties in the universal enveloping
algebra for a simple Lie algebra $\g$ and in some modules over $\g$
were constructed by Lakshmibai, Musili and Seshadri~\cite{lms:gg},
Littelmann~\cite{l:ac, l:cc}, Chari and Xi~\cite{cv:mb}
({\it monomial\/} bases); De~Concini and Kazhdan~\cite{ck:sb}, Xi~\cite{x:sb}
({\it special\/} bases and their $q$-analogs); 
Gelfand and Zelevinsky~\cite{gz:mp},
Retakh and Zelevinsky~\cite{rz:ba}, Mathieu~\cite{m:gb} ({\it good\/} bases);
Lusztig~\cite{l:cb}, Kashiwara~\cite{k:cq},
Du~\cite{d:cb1, d:cb2} ({\it canonical\/} or {\it crystal\/} bases);
see also Mathieu~\cite{m:bre} for a review and more 
references. Algorithms for computing the global crystal bases
of irreducible modules for the classical Lie algebras were recently
given by Leclerc and Toffin~\cite{lt:sa} and Lecouvey~\cite{l:aco, l:acs}.
In general, no explicit formulas are known,	however,
for the matrix elements of the generators in such
bases other than those of Gelfand and Tsetlin type.
It is known, although, that for the canonical
bases the matrix elements of the standard generators are nonnegative integers.
Some classes of representations of the symplectic, odd orthogonal and 
the Lie algebras of type $G_2$
were explicitly constructed by Donnelly~\cite{d:sa, d:ec, d:ep}
and Donnelly, Lewis and Pervine~\cite{dlp:cr}. 
The constructions were
applied to establish combinatorial properties of the
supporting graphs of the representations and were 
inspired by the earlier results of
Proctor~\cite{p:rs, p:bl, p:ss}. 
Another graph-theoretic approach is developed by 
Wildberger~\cite{w:cc, w:ccg, w:mp, w:qd} to construct simple Lie algebras
and their minuscule representations; see also Stembridge~\cite{s:mr}.

We now discuss the main idea which leads to the 
construction of the Gelfand--Tsetlin bases. 
The first point is to regard a given classical Lie algebra
not as a single object but as a part of a chain
of subalgebras with natural embeddings.
We illustrate this idea using representations of the 
symmetric groups $\Sym_n$ as an example.
Consider the chain of subgroups
\beql{chainsn}
\Sym_1\subset \Sym_2\subset\cdots\subset \Sym_n,
\end{equation}
where the subgroup $\Sym_k$ of $\Sym_{k+1}$ consists of the permutations
which fix the index $k+1$ of the set $\{1,2,\dots,k+1\}$.
The irreducible representations of the group $\Sym_n$
are indexed by partitions $\la$ of $n$. A partition 
$\la=(\la_1,\dots,\la_l)$ with
$\la_1\geq \la_2\geq\cdots\geq\la_l$
is depicted graphically as a {\it Young diagram\/}
which consists of $l$ left-justified rows of boxes so that the top row
contains $\la_1$ boxes, the second row $\la_2$ boxes, etc.
Denote by $V(\la)$ the irreducible representation of $\Sym_n$ corresponding
to the partition $\la$. One of the central results of the representation theory
of the symmetric groups is the following {\it branching rule\/}
which describes the restriction of $V(\la)$ to the subgroup $\Sym_{n-1}$:
\beq
V(\la)|^{}_{\Sym_{n-1}}=\underset{\mu}\oplus\ts	V'(\mu),
\non
\end{equation}
summed over all partitions $\mu$ whose Young diagram is obtained from that of
$\la$ by removing one box. Here
$V'(\mu)$ denotes the irreducible representation of $\Sym_{n-1}$ corresponding
to a partition $\mu$.  Thus, the restriction of $V(\la)$ to $\Sym_{n-1}$
is {\it multiplicity-free\/}, i.e., is contains
each irreducible representation	of $\Sym_{n-1}$ at most once.
This makes it possible to obtain a natural parameterization
of the basis vectors in $V(\la)$ by taking its further restrictions
to the subsequent subgroups of the chain \eqref{chainsn}. Namely, the basis
vectors will be parametrized by sequences of partitions
\beql{seqpart}
\la^{(1)}\to \la^{(2)}\to\cdots\to \la^{(n)}=\la,
\non\end{equation}
where $\la^{(k)}$ is obtained from $\la^{(k+1)}$ by removing one box.
Equivalently, each sequence of this type can be regarded as
a {\it standard tableau of shape $\la$\/} 
which is obtained by writing the numbers
$1,\dots,n$ into the boxes of $\la$ in such a way that the numbers
increase along the rows and down the columns.
In particular, the dimension of $V(\la)$ equals the number of standard
tableaux of shape $\la$.
There is only one irreducible
representation of the trivial group $\Sym_1$ therefore the procedure
defines basis vectors up to a scalar factor.
The corresponding basis is called the {\it Young basis\/}.
The symmetric group $\Sym_n$ is generated by the adjacent
transpositions $s_i=(i,i+1)$. The construction of the representation $V(\la)$
can be completed by deriving explicit formulas for the action of
the elements $s_i$ in the basis which are also due to A.~Young.
This realization of $V(\la)$ is usually called 
{\it Young's orthogonal\/} (or {\it seminormal\/}) {\it form\/}.
The details can be found, e.g., in 
James and Kerber~\cite{jk:rt} and Sagan~\cite{s:sg};
see also Okounkov and Vershik~\cite{ov:na} where an alternative
construction of the Young basis is produced.
Branching rules
and the corresponding Bratteli diagrams were employed by 
Halverson and Ram~\cite{hr:ca},	Leduc and Ram~\cite{lr:rh},
Ram~\cite{r:sr} to compute irreducible representations
of the Iwahori--Hecke algebras and some families
of centralizer algebras.

Quite a similar method can be applied to representations of the
classical Lie algebras. 
Consider the
general linear Lie algebra $\gl_n$ which consists of complex $n\times n$-matrices
with the usual matrix commutator. The chain \eqref{chainsn} is now replaced by
\beq
\gl_1\subset \gl_2\subset\cdots\subset \gl_n,
\non\end{equation}
with natural embeddings $\gl_k\subset\gl_{k+1}$.
The orthogonal Lie algebra $\oa_N$ can be regarded as a subalgebra
of $\gl_N$ which consists of skew-symmetric matrices.
Again, we have a natural chain
\beql{chainon}
\oa_2\subset \oa_3\subset\cdots\subset \oa_N.
\end{equation}
Both restrictions $\gl_n\downarrow \gl_{n-1}$ and $\oa_N\downarrow \oa_{N-1}$
are multiplicity-free so that the application of the argument which
we used for the chain \eqref{chainsn} produces basis vectors 
in an irreducible representation of $\gl_n$ or $\oa_N$. With an
appropriate normalization, these bases are precisely those
of Gelfand and Tsetlin given in \cite{gt:fdu} and \cite{gt:fdo}.
Instead of the standard tableaux, the basis vectors here are
parametrized by combinatorial objects called the {\it Gelfand--Tsetlin patterns\/}.

However, this approach does not work for the symplectic Lie algebras 
$\spa_{2n}$ since the restriction $\spa_{2n}\downarrow \spa_{2n-2}$
is not multiplicity-free. The multiplicities
are given by Zhelobenko's branching rule~\cite{z:cg} which was
re-discovered later by 
Hegerfeldt~\cite{h:bt}\footnote{Some western authors refered 
to Hegerfeldt's result as the original derivation of the rule.}.
Various attempts to fix this problem
were made by several authors. 
A natural idea is to introduce 
an intermediate Lie algebra ``$\spa_{2n-1}$"
and try to restrict an irreducible representation of $\spa_{2n}$
first to this subalgebra and then to $\spa_{2n-2}$ in the hope
to get simple spectra in the two restrictions.
Such intermediate subalgebras and their representations were studied by
Gelfand--Zelevinsky~\cite{gz:mr}, Proctor~\cite{p:os}, Shtepin~\cite{s:il}.
The drawback of this approach is the fact that the Lie algebra
$\spa_{2n-1}$ is not reductive so that
the restriction of an irreducible representation of $\spa_{2n}$
to $\spa_{2n-1}$ is not completely reducible.
In some sense, the separation of multiplicities can be achieved by constructing
a filtration of $\spa_{2n-1}$-modules; cf. Shtepin~\cite{s:il}.

Another idea is to use the restriction $\gl_{2n}\downarrow \spa_{2n}$.
Gould and Kalnins~\cite{g:rt, gk:pb} constructed
a basis for the representations of the symplectic Lie algebras
parametrized by a subset of 
the Gelfand--Tsetlin
$\gl_{2n}$-patterns. Some matrix element formulas are also derived
by using the $\gl_{2n}$-action.
A similar observation is made independently by Kirillov~\cite{k:rg}
and Proctor~\cite{p:os}. A description of the Gelfand--Tsetlin
patterns for $\spa_{2n}$ and $\oa_N$ can be obtained by
regarding them as fixed points of involutions of the
Gelfand--Tsetlin
patterns for the corresponding Lie algebra $\gl_N$.

The lowering operators in the symplectic case were given
by Mickelsson~\cite{m:lo}; see also Bincer~\cite{b:ml, b:mlo}.
The application of ordered monomials
in the lowering operators to the highest vector yields a basis
of the representation. However, the action of the Lie algebra 
generators in such a basis
does not seem to be computable. The reason is the fact that,
unlike the cases
of $\gl_n$ and $\oa_N$, the lowering operators do not
commute so that the basis depends on the chosen ordering.
A ``hidden symmetry" has been needed (cf. Cherednik~\cite{c:qg})
to make a natural choice
of an appropriate combination of the lowering operators. 
New ideas which led to a construction of a Gelfand--Tsetlin type basis
for any irreducible finite-dimensional representation of $\spa_{2n}$
came from the theory of {\it quantized enveloping algebras\/}. This is a part of
the theory of {\it quantum groups\/} originated from the works of
Drinfeld~\cite{d:ha, d:qg} and Jimbo~\cite{j:qd}.
A particular class of quantized enveloping algebras called
{\it twisted Yangians\/} introduced by Olshanski~\cite{o:ty} 
plays the role of the hidden symmetries
for the construction of the basis. We refer the reader to the book by
Chari and Pressley~\cite{cp:gq} and the review papers~\cite{m:ya} and \cite{mno:yc}
for detailed expositions of the properties of these algebras and their
origins. For each classical Lie algebra we attach the {\it Yangian\/}
$\Y(N)=\Y(\gl_N)$, or the {\it twisted Yangian\/} $\Y^{\pm}(N)$ as follows
\beq
\bal
{\rm type\ } A_n&\qquad\qquad   {\rm type}\ B_n
\qquad\qquad {\rm type}\  C_n\qquad\qquad {\rm type}\  D_n\\
\Y(n+1)&\qquad\quad	\Y^+(2n+1)\qquad\qquad   \Y^-(2n)\qquad\qquad  \Y^+(2n).\quad
\eal
\non\end{equation}
The algebra $\Y(N)$ was first introduced in the work 
of Faddeev
and the St.-Petersburg school in relation with the {\it inverse
scattering method}; see for instance
Takhtajan--Faddeev~\cite{tf:qi}, Kulish--Sklyanin~\cite{ks:qs}.
Olshanski~\cite{o:ty} introduced the twisted Yangians
in relation with his {\it centralizer construction\/}; see also \cite{mo:cc}.
In particular, he established the following key fact which plays an important role
in the basis construction. 
Given irreducible representations $V(\la)$ and $V'(\mu)$
of $\spa_{2n}$ and $\spa_{2n-2}$, respectively, there exists a natural 
irreducible action
of the algebra $\Y^-(2)$ on the space $\Hom_{\spa_{2n-2}}(V'(\mu),V(\la))$.
The homomorphism space is isomorphic to the subspace
$V(\la)^+_{\mu}$ of $V(\la)$ which is spanned
by the highest vectors of weight $\mu$ for the subalgebra $\spa_{2n-2}$.
Finite-dimensional irreducible representations
of the  twisted Yangians
were classified later in \cite{m:fd}. In particular, it turned
out that the representation $V(\la)^+_{\mu}$ of $\Y^-(2)$
can be extended to the Yangian $\Y(2)$.	Another proof of this fact was
given recently by Nazarov~\cite{n:rt}.
The algebra $\Y(2)$ and its representations are well-studied;
see Tarasov~\cite{t:im}, Chari--Pressley~\cite{cp:yr}. 
A large class of representation of $\Y(2)$ admits
Gelfand--Tsetlin-type bases associated with the
inclusion $\Y(1)\subset\Y(2)$; see \cite{m:gt}.
This allows one to get a natural basis in the space
$V(\la)^+_{\mu}$
and then by induction to get a basis in the entire space $V(\la)$.
Moreover, it turns out to be possible to write down explicit formulas
for the action of the generators of the symplectic Lie
algebra in this basis; see the author's paper \cite{m:br}
for more details. This construction together with
the work of Gelfand and Tsetlin thus provides explicit
realizations of all finite-dimensional irreducible
representations of the classical Lie algebras.

The same method can be applied to the pairs of the orthogonal Lie algebras
$\oa_{N-2}\subset\oa_{N}$. Here the corresponding space $V(\la)^+_{\mu}$
is a natural $\Y^+(2)$-module which can also be extended to
a $\Y(2)$-module. This leads to a construction of a natural basis
in the representation $V(\la)$ and allows one to explicitly
calculate the representation matrices; see \cite{m:wb, m:wbg}.
This realization of $V(\la)$ is alternative to that
of Gelfand and Tsetlin~\cite{gt:fdo}. To compare the two constructions,
note that the basis of \cite{gt:fdo}
in the orthogonal case  lacks the {\it weight\/} property, 
i.e., the basis vectors
are not eigenvectors for the Cartan subalgebra. The reason for that is the fact
that the chain \eqref{chainon} involves Lie algebras of different types
($B$ and $D$) and the embeddings are not compatible with the
root systems. In the new approach we use instead the chains
\beq
\oa_2\subset \oa_4\subset\cdots\subset \oa_{2n}\qquad\text{and}\qquad
\oa_3\subset \oa_5\subset\cdots\subset \oa_{2n+1}.
\non\end{equation}
The embeddings here ``respect" the root systems so that the basis of $V(\la)$ possesses
the weight property in both the symplectic and orthogonal cases.
However, the new weight bases, in their turn, lack the {\it orthogonality\/}
property of the Gelfand--Tsetlin bases: the latter are orthogonal
with respect to the standard inner product in the 
representation space $V(\la)$. It is an open problem to construct
a natural basis of $V(\la)$ in the $B,C$ and $D$ cases
which would simultaneously accommodate the
two properties.

This chapter is structured as follows. In Section~\ref{sec:gtgln} 
we review the construction
of the Gelfand--Tsetlin basis for the general linear Lie algebra
and discuss its various versions. We start by applying 
the most elementary approach which consists of using
explicit formulas for the lowering operators in a way similar to
the pioneering works of the sixties.
Remarkably, these operators admit several other
presentations which reflect different approaches to the problem
developed in the literature. First, we outline the general theory
of extremal projectors and Mickelsson algebras as a natural way to
work with the lowering operators. Next, we describe the $\gl_n$-type 
{\it Mickelsson--Zhelobenko algebra\/} which is then used to
prove the branching rule and derive the matrix element formulas. 
Further, we outline the Gould construction based upon the
characteristic identities. Finally, we produce quantum minor
formulas for the lowering operators inspired by the Yangian approach
and describe the action of the Drinfeld 
generators in the Gelfand--Tsetlin basis.

In Section~\ref{sec:osp} we produce
weight bases for representations of the orthogonal and symplectic Lie algebras.
Here we describe the relevant Mickelsson--Zhelobenko algebra,
formulate the branching rules and construct the basis vectors.
Then we outline
the properties of the (twisted) Yangians and their representations
and explain their relationship with the lowering and raising operators.
Finally, we sketch the main ideas in the calculation 
of the matrix element formulas.

Section~\ref{sec:gtoN} is devoted to the Gelfand--Tsetlin bases
for the orthogonal Lie algebras. We outline the basis construction
along the lines of the general method of
Mickelsson algebras.

At the end of each section we give brief bibliographical comments
pointing towards the original articles and to
the references where the proofs or further details
can be found.

It gives me pleasure to thank I.~M.~Gelfand
for his comment on the preliminary version of the paper.
My thanks also extend to
V.~K.~Dobrev,
V.~M.~Futorny,
M.~D.~Gould,
M.~Harada,
M.~L.~Nazarov, 
G.~I.~Olshanski,
S.~A.~Ovsienko,
T.~D.~Palev,
V.~S.~Retakh, and
V.~N.~Tolstoy 
who sent me remarks and references.

\section{Gelfand--Tsetlin basis for representations of $\gl_n$}\label{sec:gtgln}
\setcounter{equation}{0}

Let $E_{ij}$, $i,j=1,\dots,n$ denote the standard basis of the general
linear Lie algebra
$\gl_n$ over the field of complex numbers.
The subalgebra $\gl_{n-1}$ is spanned by
the basis elements $E_{ij}$ with $i,j=1,\dots,n-1$. Denote by $\h=\h_n$ the diagonal
Cartan subalgebra in $\gl_n$. The elements $E_{11}, \dots,E_{nn}$ form
a basis of $\h$.

Finite-dimensional irreducible representations of $\gl_n$
are in a one-to-one correspondence with $n$-tuples
of complex numbers $\lambda=(\lambda_1,\dots,\lambda_n)$ 
such that
\beql{cond}
\lambda_i-\lambda_{i+1}\in\ZZ_+\qquad\text{for}\quad i=1,\dots,n-1.
\end{equation}
Such an $n$-tuple $\lambda$ is called the {\it highest weight\/}
of the corresponding representation which
we shall denote by $L(\lambda)$.
It contains a unique, up to a multiple, nonzero vector $\xi$
(the {\it highest vector\/}) such that
$
E_{ii}\ts\xi=\lambda_i\ts\xi
$
for $i=1,\dots,n$ and
$
E_{ij}\ts\xi=0
$
for $1\leq i<j\leq n$. 

The following theorem is the {\it branching rule\/}
for the reduction $\gl_{n}\downarrow\gl_{n-1}$.

\bth\label{thm:brar}
The restriction of $L(\la)$ to the subalgebra $\gl_{n-1}$ is isomorphic
to the direct sum of pairwise inequivalent irreducible representations
\beq
L(\lambda)|^{}_{\gl_{n-1}}\simeq \underset{\mu}\oplus\ts L'(\mu),
\non
\end{equation}
summed over the highest
weights $\mu$ satisfying the betweenness conditions
\beq\label{amulam2}
\lambda_i-\mu_i\in\ZZ_+ 	\quad\text{and}\quad \mu_i-\lambda_{i+1}\in
\ZZ_+ 	\qquad\text{for}\quad i=1,\dots,n-1.
\end{equation}
\eth

The rule could presumably be attributed to I.~Schur who
was the first to discover the representation-theoretic significance
of a particular class of symmetric polynomials which now bear his name.
Without loss of generality we may regard $\la$ as a partition:
we can take the composition of $L(\la)$
with an appropriate automorphism of $\U(\gl_n)$ 
which sends $E_{ij}$ to $E_{ij}+\delta_{ij}\ts a$
for some $a\in\CC$.
The {\it character\/} of $L(\la)$ regarded as a $GL_n$-module
is the {\it Schur polynomial\/} $s_{\la}(x)$, $x=(x_1,\dots,x_n)$ defined by
\beq
s_{\la}(x)=\tr\big(g,\ts L(\la)\big),
\non\end{equation}
where $x_1,\dots,x_n$ are the eigenvalues of $g\in GL_n$.
The Schur polynomial is
symmetric in the $x_i$ and can be given
by the explicit combinatorial formula
\beql{schur}
s_{\la}(x)=\sum_{T}x^T,
\end{equation}
summed over the {\it semistandard tableaux\/} $T$ of shape $\la$
(cf. Remark~\ref{rem:part} below), where
$x^T$ is the monomial containing $x_i$ with the power equal to
the number of occurrences of $i$ in $T$; see, e.g., 
Macdonald~\cite[Chapter~1]{m:sf}
or Sagan~\cite[Chapter~4]{s:sg}
for more details. To find out what happens when $L(\la)$ is restricted
to $GL_{n-1}$ we just need to put $x_n=1$ into the formula \eqref{schur}.
The right hand side will then be written as the sum of the Schur polynomials
$s_{\mu}(x_1,\dots,x_{n-1})$ with $\mu$ satisfying \eqref{amulam2}.

On the other hand, the multiplicity-freeness of the reduction 
$\gl_n\downarrow\gl_{n-1}$ can be explained by the fact that
the vector space $\Hom_{\gl_{n-1}}(L'(\mu),L(\la))$ bears
a natural irreducible representation of the 
centralizer $\U(\gl_n)^{\gl_{n-1}}$;
see, e.g., Dixmier~\cite[Section~9.1]{d:ae}. However,
the centralizer is a commutative algebra and therefore
if the homomorphism space is nonzero then it must be one-dimensional.

The branching rule is implicit in the 
formulas of Gelfand and Tsetlin~\cite{gt:fdu}. Its proof based
upon an explicit realization of the representations of $GL_n$
was given by Zhelobenko~\cite{z:cg}. We outline
a proof of Theorem~\ref{thm:brar} below in Section~\ref{subsec:ama}
which employs the modern theory of Mickelsson algebras
following Zhelobenko~\cite{z:rr}. Two other proofs can be found
in Goodman and Wallach~\cite[Chapters~8 \& 12]{gw:ri}.

The subsequent applications of the branching rule to the subalgebras
of the chain
\beq
\gl_1\subset\gl_2\subset\cdots\subset\gl_{n-1}\subset\gl_n
\non\end{equation}
yield a parameterization of basis vectors in $L(\lambda)$ by
the combinatorial objects called
the {\it Gelfand--Tsetlin patterns\/}. 
Such a pattern $\Lambda$ (associated with
$\lambda$) is an array of row vectors
\begin{align}
&\qquad\lambda^{}_{n1}\qquad\lambda^{}_{n2}
\qquad\qquad\cdots\qquad\qquad\lambda^{}_{nn}\non\\
&\qquad\qquad\lambda^{}_{n-1,1}\qquad\ \ \cdots\ \ 
\ \ \qquad\lambda^{}_{n-1,n-1}\non\\
&\quad\qquad\qquad\cdots\qquad\cdots\qquad\cdots\non\\
&\quad\qquad\qquad\qquad\lambda^{}_{21}\qquad\lambda^{}_{22}\non\\
&\quad\qquad\qquad\qquad\qquad\lambda^{}_{11}  \non
\end{align}
where the upper row coincides with $\lambda$ and 
the following conditions hold
\beq\label{aconl}
\lambda^{}_{ki}-\lambda^{}_{k-1,i}\in\ZZ_+,\qquad  
\lambda^{}_{k-1,i}-\lambda^{}_{k,i+1}\in\ZZ_+,\qquad 
i=1,\dots,k-1
\end{equation}
for each $k=2,\dots,n$. 

\bre\label{rem:part}
{\rm
If the highest weight $\lambda$ is a partition then
there is a natural
bijection between the patterns associated with $\lambda$ and
the {\it semistandard\/} $\lambda$-tableaux with entries in $\{1,\dots,n\}$.
Namely, the pattern $\La$ can be viewed as the sequence of partitions
\beq
\lambda^{(1)}\subseteq \lambda^{(2)}\subseteq\cdots 
\subseteq \lambda^{(n)}=\lambda, \non
\end{equation}
with 
$\lambda^{(k)}=(\lambda^{}_{k1},\dots,\lambda^{}_{kk})$. 
Conditions \eqref{aconl} mean that the skew diagram
$\lambda^{(k)}/\lambda^{(k-1)}$ is a {\it horizontal strip\/};
see, e.g., Macdonald~\cite[Chapter~1]{m:sf}.
The corresponding semistandard tableau is obtained by placing
the entry $k$ into each box of $\lambda^{(k)}/\lambda^{(k-1)}$.
}
\ere


The {\it Gelfand--Tsetlin basis\/} of $L(\la)$ is 
provided by the following theorem.
Let us set
$l^{}_{ki}=\lambda^{}_{ki}-i+1$.

\bth\label{thm:abasis} There exists a basis 
$\{\xi^{}_{\Lambda}\}$ in $L(\lambda)$ parametrized by all
patterns $\Lambda$ such that the action
of generators of $\gl_n$ is given by the formulas
\begin{align}\label{aekk}
E_{kk}\ts \xi^{}_{\Lambda}&=\left(\sum_{i=1}^k\lambda^{}_{ki}
-\sum_{i=1}^{k-1}\lambda^{}_{k-1,i}\right)
\xi^{}_{\Lambda},
\\
\label{aekk+1}
E_{k,k+1}\ts \xi^{}_{\Lambda}
&=-\sum_{i=1}^k \frac{(l^{}_{ki}-l^{}_{k+1,1})\cdots (l^{}_{ki}-l^{}_{k+1,k+1})}
{(l^{}_{ki}-l^{}_{k1})\cdots \wedge\cdots(l^{}_{ki}-l^{}_{kk})}
\ts
\xi^{}_{\Lambda+\delta^{}_{ki}},\\
\label{aek+1k}
E_{k+1,k}\ts \xi^{}_{\Lambda}
&=\sum_{i=1}^k \frac{(l^{}_{ki}-l^{}_{k-1,1})\cdots (l^{}_{ki}-l^{}_{k-1,k-1})}
{(l^{}_{ki}-l^{}_{k1})\cdots \wedge\cdots(l^{}_{ki}-l^{}_{kk})}
\ts
\xi^{}_{\Lambda-\delta^{}_{ki}}.
\end{align}
The arrays $\Lambda\pm\delta^{}_{ki}$
are obtained from $\Lambda$ by replacing $\lambda^{}_{ki}$
by $\lambda^{}_{ki}\pm1$. It is supposed
that $\xi^{}_{\Lambda}=0$ if the array $\Lambda$ is not a pattern;
the symbol $\wedge$ indicates that the zero factor in the denominator
is skipped.
\eth

A construction of the basis vectors is given in Theorem~\ref{thm:abasislo}
below. A derivation of the matrix element formulas \eqref{aekk}--\eqref{aek+1k}
is outlined in Section~\ref{subsec:ama}.

The vector space $L(\la)$ is equipped with a contravariant
inner product $\lan,\ran$.
It is uniquely 
determined by the conditions
\beq
\lan \xi,\xi \ran =1\qquad\text{and}\qquad
\lan E_{ij}\ts\eta, \zeta \ran =  \lan \eta, E_{ji}\ts\zeta \ran
\non\end{equation}
for any vectors $\eta,\zeta\in L(\la)$ and any indices $i,j$.
In other words, for the adjoint operator for $E_{ij}$
with respect to the inner product we have $(E_{ij})^*=E_{ji}$. 

\bpr\label{prop:norms}
The basis $\{\xi^{}_{\Lambda}\}$ is orthogonal with respect to the
inner product $\lan,\ran$.
Moreover, we have
\beq
\lan \xi^{}_{\Lambda}, \xi^{}_{\Lambda}\ran=
\prod_{k=2}^n\ts\ts\prod_{1\leq i\leq j<k}
\frac{(l_{ki}-l_{k-1,j})!}
{(l_{k-1,i}-l_{k-1,j})!}
\prod_{1\leq i<j\leq k}
\frac{(l_{ki}-l_{kj}-1)!}
{(l_{k-1,i}-l_{kj}-1)!}.
\non\end{equation}
\epr

The formulas of Theorem~\ref{thm:abasis}
can therefore be rewritten in the orthonormal basis 
\beql{orthonorm}
\zeta^{}_{\Lambda}=\xi^{}_{\Lambda}/\|\ts\xi^{}_{\Lambda}\|,\qquad
\|\ts\xi^{}_{\Lambda}\|^2= \lan \xi^{}_{\Lambda}, \xi^{}_{\Lambda}\ran.
\end{equation}
They were presented in this form
in the original work by Gelfand and Tsetlin~\cite{gt:fdu}. 
A proof of Proposition~\ref{prop:norms} will 
be outlined in Section~\ref{subsec:ama}.

\subsection{Construction of the basis: lowering and raising operators}\label{subsec:loexp}

For each $i=1,\dots,n-1$
introduce the following elements of the universal enveloping algebra
$\U(\gl_n)$ 
\begin{align}\label{zin}
z_{in}&=\sum_{i>i_1>\cdots>i_s\geq 1}
E_{ii_1}E_{i_1i_2}\cdots E_{i_{s-1}i_s}E_{i_sn}
(h_i-h_{j_1})\cdots (h_i-h_{j_r}),
\\
\label{zni}
z_{ni}&=\sum_{i<i_1<\cdots<i_s<n}
E_{i_1i}E_{i_2i_1}\cdots E_{i_si_{s-1}}E_{ni_s}
(h_i-h_{j_1})\cdots (h_i-h_{j_r}),
\end{align}
where $s$ runs over nonnegative integers, $h_i=E_{ii}-i+1$ 
and $\{j_1,\dots,j_r\}$
is the complementary subset to $\{i_1,\dots,i_s\}$
in the set $\{1,\dots,i-1\}$ or $\{i+1,\dots,n-1\}$, respectively.
For instance,
\beq
\bal
z_{13}&=E_{13},\qquad z_{23}=E_{23}\ts(h_2-h_1)+E_{21}E_{13},\\
z_{32}&=E_{32},\qquad z_{31}=E_{31}\ts(h_1-h_2)+E_{21}E_{32}.
\eal
\non\end{equation}

Consider now the irreducible finite-dimensional representation $L(\la)$
of $\gl_n$ with the highest weight $\la=(\la_1,\dots,\la_n)$ and the
highest vector $\xi$.
Denote by $L(\lambda)^+$ the subspace of $\gl_{n-1}$-highest vectors
in $L(\lambda)$:
\beq
L(\lambda)^+=\{\eta\in L(\lambda)\ |\ E_{ij}\ts \eta=0,
\qquad 1\leq i<j<n\}.
\non\end{equation}
Given a $\gl_{n-1}$-weight
$\mu=(\mu_1,\dots,\mu_{n-1})$ we denote by $L(\lambda)^+_{\mu}$
the corresponding weight subspace in $L(\lambda)^+$:
\beq
L(\lambda)^+_{\mu}=\{\eta\in L(\lambda)^+\ |\ E_{ii}\ts\eta=
\mu_i\ts\eta,\qquad i=1,\dots,n-1\}.
\non\end{equation}

The main property of the elements $z_{ni}$ and $z_{in}$ is described by
the following lemma.

\ble\label{lem:lowo}
Let $\eta\in L(\lambda)^+_{\mu}$.
Then for any $i=1,\dots,n-1$ we have
\beq
z_{in}\ts \eta\in L(\lambda)^+_{\mu+\delta_i}\qquad\text{and}
\qquad z_{ni}\ts \eta \in L(\lambda)^+_{\mu-\delta_i},
\non\end{equation}
where the weight $\mu\pm\delta_i$ is obtained from $\mu$ by replacing $\mu_i$ with
$\mu_i\pm 1$.
\ele

This result allows us to regard the elements $z_{in}$ and $z_{ni}$ as operators
in the space $L(\lambda)^+$. They are called the {\it raising\/} and {\it lowering
operators\/}, respectively. By the branching rule (Theorem~\ref{thm:brar})
the space $L(\la)^+_{\mu}$ is one-dimensional if the conditions \eqref{amulam2}
hold and it is zero otherwise.	The following lemma will be proved in 
Section~\ref{subsec:ama}.

\ble\label{lem:basms} 
Suppose that $\mu$ satisfies
the betweenness conditions \eqref{amulam2}.
Then the vector
\beq
\xi_{\mu}=z_{n1}^{\la_1-\mu_1}\cdots z_{n,n-1}^{\la_{n-1}-\mu_{n-1}}\ts\xi
\non\end{equation}
is nonzero. Moreover, the space $L(\la)^+_{\mu}$ is spanned by $\xi_{\mu}$.
\ele

The $\U(\gl_{n-1})$-span of each nonzero vector $\xi_{\mu}$
is a $\gl_{n-1}$-module isomorphic to $L'(\mu)$.
Iterating the construction of the vectors $\xi_{\mu}$
for each pair of Lie algebras $\gl_{k-1}\subset\gl_k$
we shall be able to get a basis in the entire space $L(\la)$.

\bth\label{thm:abasislo}
The basis vectors $\xi^{}_{\Lambda}$ of Theorem~\ref{thm:abasis}
can be given by the formula
\beql{axiL}
\xi^{}_{\Lambda}=\prod_{k=2,\dots,n}^{\rightarrow}
\Bigl(z_{k1}^{\lambda^{}_{k1}-\lambda^{}_{k-1,1}}\cdots
z_{k,k-1}^{\lambda^{}_{k,k-1}-\lambda^{}_{k-1,k-1}}\Bigr)\ts\xi,
\end{equation}
where the factors in the product are ordered in accordance with
increase of the indices.
\eth

\subsection{The Mickelsson algebra theory}\label{subsec:ma}

The lowering and raising
operators $z_{ni}$ and $z_{in}$ in the space $L(\lambda)^+$ (see Lemma~\ref{lem:lowo})
satisfy some quadratic relations with rational coefficients
in the parameters of the highest weights. These relations
can be regarded in a representation independent form with a suitable
interpretation of the coefficients as rational functions 
in the elements of the Cartan subalgebra $\h$. 
In such an abstract form
the algebras of lowering and raising
operators were introduced by Mickelsson~\cite{m:sa} who, however, did not
use any rational extensions of the algebra $\U(\h)$. The importance
of this extension was realized by Zhelobenko~\cite{z:sa, z:za} who developed
a general structure theory of these algebras which he called
the {\it Mickelsson algebras\/}. Another important ingredient is the theory
of {\it extremal projectors\/} originated from the works
of Asherova, Smirnov and Tolstoy \cite{ast:po, ast:po2, ast:dc} and
further developed by Zhelobenko~\cite{z:it, z:rr}.

Let $\g$ be a Lie algebra over $\CC$ and $\kgot$ be its subalgebra
reductive in $\g$. This means that the adjoint $\kgot$-module $\g$
is completely reducible. In particular, $\kgot$ is a reductive
Lie algebra. Fix a Cartan subalgebra $\h$ of $\kgot$ and 
a triangular
decomposition
\beq
\kgot=\kgot^-\oplus\h\oplus\kgot^+.
\non\end{equation}
The subalgebras $\kgot^-$ and $\kgot^+$ are respectively spanned by the
negative and positive root vectors $e_{-\al}$ and $e_{\al}$ with
$\al$ running over the set of positive roots $\Delta^+$
of $\kgot$ with respect to $\h$. The root vectors will be assumed
to be normalized in such a way that
\beql{normalirv}
[e_{\al},e_{-\al}]=h_{\al},\qquad \al(h_{\al})=2
\end{equation}
for all $\al\in\Delta^+$.

Let $\J=\U(\g)\ts\kgot^+$ 
be the left ideal of $\U(\g)$
generated by $\kgot^+$. Its normalizer $\Norm\J$ is a subalgebra
of $\U(\g)$ defined by
\beq
\Norm\J=\{u\in\U(\g)\ |\ \J\ts u\subseteq \J\}.
\non\end{equation}
Then $\J$ is a two-sided ideal of $\Norm\J$ and the {\it Mickelsson algebra\/}
$\Sr(\g,\kgot)$ is defined as the quotient
\beq
\Sr(\g,\kgot)=\Norm\J/\J.
\non\end{equation}

Let $\R(\h)$ denote
the field of fractions of the commutative algebra $\U(\h)$.
In what follows
it is convenient to consider the extension $\U'(\g)$ of the
universal enveloping algebra $\U(\g)$ defined by
\beq
\U'(\g)=\U(\g)\ot_{\U(\h)} \R(\h).
\non\end{equation}
Let $\J'=\U'(\g)\ts\kgot^+$ 
be the left ideal of $\U'(\g)$
generated by $\kgot^+$. Exactly as with the ideal $\J$ above,
$\J'$ is a two-sided ideal of the normalizer $\Norm\J'$ and the 
{\it Mickelsson--Zhelobenko algebra\/}\footnote{Zhelobenko sometimes used the names
{\it Z-algebra\/} or {\it extended Mickelsson algebra\/}.
The author believes the new name is more emphatic
and justified from the scientific point of view.}
$\Z(\g,\kgot)$ is defined as the quotient
\beq
\Z(\g,\kgot)=\Norm\J'/\J'.
\non\end{equation}
Clearly, $\Z(\g,\kgot)$ is an extension of the Mickelsson algebra $\Sr(\g,\kgot)$,
\beq
\Z(\g,\kgot)=\Sr(\g,\kgot)\ot_{\U(\h)} \R(\h).
\non\end{equation}
An equivalent definition of the algebra $\Z(\g,\kgot)$ can be given by using
the quotient space
\beq
\M(\g,\kgot)= \U'(\g)/\J'.
\non\end{equation}
The Mickelsson--Zhelobenko algebra $\Z(\g,\kgot)$ coincides with the subspace
of $\kgot$-highest vectors in $\M(\g,\kgot)$
\beq
\Z(\g,\kgot)=\M(\g,\kgot)^+,
\non\end{equation}
where
\beq
\M(\g,\kgot)^+=\{v\in \M(\g,\kgot)\ |\ \kgot^+\ts v=0\}.
\non\end{equation}

The algebraic structure of the algebra $\Z(\g,\kgot)$ can be described with
the use of the {\it extremal projector\/} for the Lie algebra $\kgot$.
In order to define it, suppose that the positive roots
are $\Delta^+=\{\al_1,\dots,\al_m\}$. Consider the vector space 
$\F_{\mu}(\kgot)$ of formal series
of weight $\mu$ monomials
\beq
e_{-\al_1}^{k_1}\cdots e_{-\al_m}^{k_m}
e_{\al_m}^{r_m}\cdots e_{\al_1}^{r_1}
\non\end{equation}
with coefficients in $\R(\h)$,
where
\beq
(r_1-k_1)\ts\al_1+\cdots+(r_m-k_m)\ts\al_m=\mu.
\non\end{equation}
Introduce the space $\F(\kgot)$ as the direct sum
\beq
\F(\kgot)=\underset{\mu}{\oplus}\ts \F_{\mu}(\kgot).
\non\end{equation}
That is, the elements of $\F(\kgot)$ are finite sums $\sum x_{\mu}$ with 
$x_{\mu}\in \F_{\mu}(\kgot)$.
It can be shown that  $\F(\kgot)$
is an algebra with respect to the natural multiplication of formal series.
The algebra $\F(\kgot)$ is equipped with a Hermitian
anti-involution (antilinear involutive anti-automorphism)
defined by
\beq
e_{\al}^*=e_{-\al},\qquad \al\in\Delta^+.					 
\non\end{equation}
Further, call an ordering of the positive roots {\it normal\/}
if any composite root lies between its components. For instance,
there are precisely two normal orderings for the root system of type $B_2$,
\beq
\Delta^+=\{\al,\al+\be,\al+2\be,\be\}\qquad\text{and}\qquad
\Delta^+=\{\be,\al+2\be,\al+\be,\al\},
\non\end{equation}
where $\al$ and $\be$ are the simple roots. In general, the number
of normal orderings coincides with the number of reduced decompositions
of the longest element of the corresponding Weyl group.

For any $\al\in \Delta^+$ introduce the element of $\F(\kgot)$ by
\beql{palpha}
p_{\al}=1+\sum_{k=1}^{\infty} e_{-\al}^k\ts e_{\al}^k
\ts\frac{(-1)^k}{k!\ts(h_{\al}+\rho(h_{\al})+1)\cdots (h_{\al}+\rho(h_{\al})+k)},
\end{equation}
where $h_{\al}$ is defined in \eqref{normalirv} and $\rho$ is the half sum
of the positive roots.	Finally, define the {\it extremal projector\/}
$p=p_{\ts\kgot}$ by
\beq
p= p_{\al_1}\cdots p_{\al_m}
\non\end{equation}
with the product taken in a normal ordering of the positive roots $\al_i$.

\bth\label{thm:extrpr}
The element $p\in \F(\kgot)$ does not depend on the normal ordering
on $\Delta^+$ and satisfies the conditions
\beql{extprpr}
e_{\al}\ts p=p\ts e_{-\al}=0\qquad\text{for all}\qquad \al\in\Delta^+.
\end{equation}
Moreover, $p^*=p$ and $p^2=p$.
\eth

In fact, the relations \eqref{extprpr} uniquely determine
the element $p$, up to a factor from $\R(\h)$. 
The extremal projector naturally acts on the vector space $\M(\g,\kgot)$.
The following corollary states that the Mickelsson--Zhelobenko algebra coincides
with its image.

\bco\label{cor:maim}
We have
\beq
\Z(\g,\kgot)=p\ts \M(\g,\kgot).
\non\end{equation}
\eco

To get a more precise description of the algebra $\Z(\g,\kgot)$ consider
a $\kgot$-module decomposition
\beq
\g=\kgot\oplus\pgot.
\non\end{equation}
Choose a weight basis $e_1,\dots,e_n$ 
(with respect to the adjoint action of $\h$)
of the complementary
module $\pgot$.

\bth\label{thm:pbwma}
The elements 
\beq
a_i=p\ts e_i,\qquad i=1,\dots,n
\non\end{equation}
are generators of the Mickelsson--Zhelobenko algebra $\Z(\g,\kgot)$.
Moreover, the monomials
\beq
a_1^{k_1}\cdots a_n^{k_n},\qquad k_i\in\ZZ_+,
\non\end{equation}
form a basis of $\Z(\g,\kgot)$.
\eth

It can be proved that the generators $a_i$ of $\Z(\g,\kgot)$
satisfy quadratic defining relations; see \cite{z:it}.
For the pairs $(\g,\kgot)$	relevant to the constructions
of bases of the Gelfand--Tsetlin type,
the relations can be explicitly written down; 
cf. Sections~\ref{subsec:ama} and \ref{subsec:rl} below. 

Regarding  $\Z(\g,\kgot)$ as a right $\R(\h)$-module,
it is possible to introduce the normalized elements
\beq
z_i=a_i\ts\pi_i,\qquad \pi_i\in\U(\h)
\non\end{equation}
by multiplying $a_i$ by its {\it right denominator\/} $\pi_i$.
Therefore the $z_i$ can be viewed as elements of the Mickelsson algebra
$\Sr(\g,\kgot)$. 

To formulate the final theorem of this section, for any $\g$-module $V$ set
\beq
V^+=\{v\in V\ |\ \kgot^+\ts v=0\}.
\non\end{equation}

\bth\label{thm:mick}
Let $V=\U(\g)\ts v$ be a cyclic $\U(\g)$-module generated by
an element $v\in V^+$. Then the subspace $V^+$
is linearly spanned by the elements
\beq
z_1^{k_1}\cdots z_n^{k_n}\ts v,\qquad k_i\in\ZZ_+.
\non\end{equation}
\eth

\subsection{Mickelsson--Zhelobenko algebra $\Z(\gl_n,\gl_{n-1})$}\label{subsec:ama}

For any positive integer $m$ consider the general linear Lie
algebra $\gl_m$. The positive roots of 
$\gl_m$ with respect to the diagonal
Cartan subalgebra $\h$ (with the standard choice
of the positive root system)
are naturally enumerated by the pairs $(i,j)$
with $1\leq i<j\leq m$.
In accordance with the general theory outlined in the previous section,
for each pair introduce the formal series $p_{ij}\in\F(\gl_m)$ by
\beq
p_{ij}=1+\sum_{k=1}^{\infty}(E_{ji})^k(E_{ij})^k\ts
\frac{(-1)^k}{k!\ts (h_i-h_j+1)\cdots (h_i-h_j+k)},
\non\end{equation}
where, as before, $h_i=E_{ii}-i+1$. Then define
the element $p=p_m$ by
\beq
p=\prod_{i<j}\ts p_{ij},
\non\end{equation}
where the product is taken in a  normal ordering on the pairs
$(i,j)$. By Theorem~\ref{thm:extrpr},
\beql{annih}
E_{ij}\ts p=p E_{ji}=0\qquad{\rm for}\quad 1\leq i<j\leq m.
\end{equation}

Now set $m=n-1$. By Theorem~\ref{thm:pbwma}, ordered
monomials in the elements $E_{nn}$,
$pE_{in}$ and $pE_{ni}$ with $i=1,\dots,n-1$  
form a basis of
$\Z(\gl_n,\gl_{n-1})$ as a left or right $\R(\h)$-module.
These elements can explicitly be given by
\begin{align}
pE_{in}&=\sum_{i>i_1>\cdots>i_s\geq 1}
E_{ii_1}E_{i_1i_2}\cdots E_{i_{s-1}i_s}E_{i_sn}
\frac{1}{(h_i-h_{i_1})\cdots (h_i-h_{i_s})},
\non\\
\label{apni}
pE_{ni}&=\sum_{i<i_1<\cdots<i_s<n}
E_{i_1i}E_{i_2i_1}\cdots E_{i_si_{s-1}}E_{ni_s}
\frac{1}{(h_i-h_{i_1})\cdots (h_i-h_{i_s})},
\end{align}
where $s=0,1,\dots$. Indeed, by choosing appropriate
normal orderings on the positive roots, we can write
\beq
pE_{in}=p_{1i}\cdots p_{i-1,i}\ts E_{in}
\qquad
\text{and} 
\qquad
pE_{ni}= p_{i,i+1}\cdots p_{i,n-1}\ts E_{ni}.
\non\end{equation}

The lowering and raising operators introduced in Section~\ref{subsec:loexp}
coincide with the normalized generators:
\begin{align}
z_{in}&=pE_{in}\ts
(h_i-h_{i-1})\cdots (h_i-h_{1}),
\non\\
\label{alow2}
z_{ni}&=pE_{ni}\ts
(h_i-h_{i+1})\cdots (h_i-h_{n-1}),
\end{align}
which belong to the Mickelsson algebra $\Sr(\gl_n,\gl_{n-1})$.
Thus, Lemma~\ref{lem:lowo} is an immediate corollary
of \eqref{annih}.

\bpr\label{prop:relama}
The lowering and raising operators satisfy the following relations
\begin{align}\label{acom}
z_{ni}z_{nj}&=z_{nj}z_{ni}\qquad\text{for all}\quad i,j,\\
\label{acom2}
z_{in}z_{nj}&=z_{nj}z_{in}\qquad\text{for}\quad i\ne j,
\end{align}
and
\beql{z_ins_ni}
z_{in}z_{ni}=\prod_{j=1,\ts j\ne i}^n(h_i-h_j-1)
+\sum_{j=1}^{n-1} z_{nj}z_{jn}
\prod_{k=1,\ts k\ne j}^{n-1} \frac{h_i-h_k-1}{h_j-h_k}.
\end{equation}
\epr

\Proof
We use the properties of $p$.
Assume that $i<j$. Then \eqref{annih} and
\eqref{apni} imply that in $\Z(\gl_n,\gl_{n-1})$
\beq
pE_{ni}\ts pE_{nj}=pE_{ni} E_{nj}, \qquad 
pE_{nj}\ts pE_{ni}=pE_{ni} E_{nj}\ts \frac{h_i-h_j+1}{h_i-h_j}.
\non\end{equation}
Now \eqref{acom} follows from \eqref{alow2}. The proof of
\eqref{acom2} is similar. The ``long" relation \eqref{z_ins_ni}
can be verified by analogous but more
complicated direct calculations. We give its different proof based upon 
the properties of the {\it Capelli determinant\/} ${\cal C}(u)$. 
Consider the $n\times n$-matrix $E$ whose
$ij$-th entry is $E_{ij}$ and let $u$ be a formal variable.
Then ${\cal C}(u)$ is a polynomial
with coefficients in the universal enveloping algebra $\U(\gl_n)$
defined by
\beql{aqdet}
{\cal C}(u)=\sum_{\sigma\in \Sym_n}\sgn \sigma\cdot (u+E)_{\sigma(1),1}\cdots
(u+E-n+1)_{\sigma(n),n}.
\end{equation}
It is well known that all its coefficients belong
to the center of $\U(\gl_n)$ and generate the center; see 
e.g. Howe--Umeda~\cite{hu:ci}.
This also easily follows from the properties of the {\it quantum
determinant\/} of the Yangian for the Lie algebra $\gl_n$;
see e.g. \cite{mno:yc}. 
Therefore, these coefficients act in $L(\lambda)$ as scalars
which can be easily found by applying ${\cal C}(u)$ to the highest vector $\xi$:
\beql{atu}
{\cal C}(u)|^{}_{L(\lambda)}=(u+l_1)\cdots (u+l_n),\qquad l_i=\la_i-i+1.
\end{equation}
On the other hand, the center of $\U(\gl_n)$ is
a subalgebra in the normalizer $\Norm \J$. 
We shall keep the same notation for the image of ${\cal C}(u)$ in the
Mickelsson--Zhelobenko algebra $\Z(\gl_n,\gl_{n-1})$.
To get explicit expressions of the coefficients of ${\cal C}(u)$ in terms of
the lowering and raising operators we
consider ${\cal C}(u)$ modulo the ideal $\J'$ and apply the projection $p$.
A straightforward calculation yields two alternative formulas
\beql{atzz}
{\cal C}(u)=(u+E_{nn})\prod_{i=1}^{n-1}(u+h_i-1)-
\sum_{i=1}^{n-1}z_{in}z_{ni}\prod_{j=1,\ts j\ne i}^{n-1}
\frac{u+h_j-1}{h_i-h_j}
\end{equation}
and
\beql{atzz2}
{\cal C}(u)=\prod_{i=1}^{n}(u+h_i)-
\sum_{i=1}^{n-1}z_{ni}z_{in}\prod_{j=1,\ts j\ne i}^{n-1}
\frac{u+h_j}{h_i-h_j}.
\end{equation}
The formulas show that ${\cal C}(u)$ can be regarded as an interpolation polynomial
for the products $z_{in}z_{ni}$ and $z_{ni}z_{in}$. Namely, 
for $i=1,\dots,n-1$ we have
\beql{capeva}
{\cal C}(-h_i+1)=(-1)^{n-1} z_{in}z_{ni}\qquad\text{and}
\qquad	{\cal C}(-h_i)=(-1)^{n-1} z_{ni}z_{in}
\end{equation}
with the agreement that when we evaluate $u$ in $\U(\h)$
we write the coefficients
of the polynomial to the left from powers of $u$.
Comparing the values of \eqref{atzz} and \eqref{atzz2} at $u=-h_i+1$
we get \eqref{z_ins_ni}. \endproof

Note that the relation inverse to \eqref{z_ins_ni} can be obtained
by comparing the values of \eqref{atzz} and \eqref{atzz2} at $u=-h_i$.

Next we outline the proofs of the branching rule
(Theorem~\ref{thm:brar}) and the formulas for the basis elements
of $L(\la)^+$ (Lemma~\ref{lem:basms}). The module $L(\la)$
is generated by the highest vector $\xi$ and we have
\beq
z_{in}\ts\xi=0,\qquad i=1,\dots,n-1.
\non\end{equation}
So, by Theorem~\ref{thm:mick}, the vector space $L(\la)^+$
is spanned by the elements
\beql{spelemlla}
z_{n1}^{k_1}\cdots z_{n,n-1}^{k_{n-1}}\ts \xi,\qquad k_i\in\ZZ_+.
\end{equation}
Let us set $\mu_i=\la_i-k_i$ for $1\leq i\leq n-1$ and denote the vector
\eqref{spelemlla} by $\xi_{\mu}$. That is,
\beql{spelemlmu}
\xi_{\mu}=z_{n1}^{\la_1-\mu_1}\cdots z_{n,n-1}^{\la_{n-1}-\mu_{n-1}}\ts 
\xi.
\end{equation}
It is now sufficient to show that the vector $\xi_{\mu}$ is nonzero
if and only if the betweenness conditions \eqref{amulam2} hold.
The linear independence of the vectors $\xi_{\mu}$
will follow from the fact that their
weights are distinct.
If $\xi_{\mu}\ne 0$ then using the relations \eqref{acom}
we conclude that each vector $z_{ni}^{\la_i-\mu_i}\ts\xi$ is nonzero.
On the other hand, $z_{ni}^{k_i}\ts\xi$ is a $\gl_{n-1}$-highest vector
of the weight obtained from $(\la_1,\dots,\la_{n-1})$ by replacing
$\la_i$ with $\la_i-k_i$. Therefore, if $k_i\geq\la_i-\la_{i+1}+1$
then the conditions \eqref{cond} are violated for this weight which implies
$z_{ni}^{k_i}\ts\xi=0$. Hence, $\la_i-\mu_i\leq \la_i-\la_{i+1}$
for each $i$, and $\mu$ satisfies \eqref{amulam2}.

For the proof of the converse statement 
we shall employ the following key lemma
which will also be used for the proof of Theorem~\ref{thm:abasis}.

\ble\label{lem:aximu}
We have for each $i=1,\dots,n-1$
\beq\label{aact}
z_{in}\ts \xi_{\mu}=-(m_i-l_1)\cdots (m_i-l_n)\ts \xi_{\mu+\delta^{}_i},
\end{equation}
where
\beq
m_i=\mu_i-i+1,\qquad l_i=\lambda_i-i+1.
\non\end{equation}
It is supposed that $\xi_{\mu+\delta^{}_i}=0$ if $\lambda_i=\mu_i$.
\ele

\Proof
The relation \eqref{acom2} implies that 
if $\la_i=\mu_i$ then $z_{in}\ts \xi_{\mu}=0$ which agrees with \eqref{aact}.
Now let $\lambda_i-\mu_i\geq 1$.
Using \eqref{acom} and \eqref{capeva}
we obtain
\beq
z_{in}\ts \xi_{\mu}=z_{in}z_{ni}\ts \xi_{\mu+\delta^{}_i}
=(-1)^{n-1}{\cal C}(-h_i+1)\ts \xi_{\mu+\delta^{}_i}
= (-1)^{n-1} {\cal C}(-m_i)\ts \xi_{\mu+\delta^{}_i}.
\non
\end{equation}
The relation \eqref{aact} now follows from \eqref{atu} and the centrality of
${\cal C}(u)$.
\endproof

If the betweenness conditions \eqref{amulam2} hold then
by Lemma~\ref{lem:aximu}, applying appropriate operators $z_{in}$ repeatedly to 
the vector $\xi_{\mu}$
we can obtain the highest vector $\xi$ with a nonzero coefficient.
This gives $\xi_{\mu}\ne 0$.

Thus, we have proved that the vectors $\xi^{}_{\La}$ defined in \eqref{axiL}
form a basis of the representation $L(\la)$. The orthogonality of the basis
vectors (Proposition~\ref{prop:norms}) is implied by the fact that
the operators $pE_{ni}$ and $pE_{in}$ are adjoint to each other
with respect to the restriction of the inner product $\lan,\ran$
to the subspace $L(\la)^+$. Therefore, for the adjoint operator to $z_{ni}$
we have
\beq
z_{ni}^*=z_{in}\frac{(h_i-h_{i+1}-1)\cdots (h_i-h_{n-1}-1)}
{(h_i-h_{1})\cdots (h_i-h_{i-1})}
\non\end{equation}
and Proposition~\ref{prop:norms} is deduced from Lemma~\ref{lem:aximu}
by induction.

Now we outline
a derivation of formulas \eqref{aekk}--\eqref{aek+1k}. 
First, since $E_{nn}\ts z_{ni}=z_{ni}\ts (E_{nn}+1)$ for any $i$,
we have
\beq
E_{nn}\ts\xi_{\mu}=\Bigl(\sum_{i=1}^n\lambda_i-
\sum_{i=1}^{n-1}\mu_i\Bigr)\ts\xi_{\mu}
\non\end{equation}
which implies \eqref{aekk}. To prove \eqref{aekk+1} is suffices
to calculate $E_{n-1,n}\ts \xi_{\mu\nu}$, where
\beq
\xi_{\mu\nu}=
z_{n-1,1}^{\mu^{}_1-\nu^{}_1}\cdots z_{n-1,n-2}^{\mu^{}_{n-2}-\nu^{}_{n-2}}
\ts\xi_{\mu}
\non\end{equation}
and the $\nu_i$ satisfy the betweenness conditions
\beq
\mu_i-\nu_i\in\ZZ_+ 	\quad\text{and}\quad \nu_i-\mu_{i+1}\in
\ZZ_+ 	\qquad\text{for}\quad i=1,\dots,n-2.
\non\end{equation}
Since $E_{n-1,n}$ commutes with the $z_{n-1,i}$, 
\beq
E_{n-1,n}\ts \xi_{\mu\nu}=
z_{n-1,1}^{\mu^{}_1-\nu^{}_1}\cdots z_{n-1,n-2}^{\mu^{}_{n-2}-\nu^{}_{n-2}}\ts
E_{n-1,n}\ts\xi_{\mu}. 
\non\end{equation}
The following lemma is implied by the explicit 
formulas for the lowering and raising operators \eqref{zin} and \eqref{zni}.

\ble\label{lem:elow}
We have the
relation in $\U'(\gl_n)$ modulo the ideal $\J'$,
\beq
E_{n-1,n}=\sum_{i=1}^{n-1}z_{n-1,i}\ts z_{in}
\frac{1}{(h_i-h_1)\cdots\wedge_i\cdots(h_i-h_{n-1})},
\non\end{equation}
where $z_{n-1,n-1}=1$.
\ele
By Lemmas~\ref{lem:aximu} and \ref{lem:elow},
\beq\label{aen-1n3}
E_{n-1,n}\ts \xi_{\mu\nu}=
-\sum_{i=1}^{n-1} \frac{(m_i-l_1)\cdots (m_i-l_n)}
{(m_i-m_1)\cdots\wedge_i\cdots (m_i-m_{n-1})}
\ts
\xi_{\mu+\delta_i,\nu}
\end{equation}
which proves \eqref{aekk+1}.
To prove \eqref{aek+1k} we use Proposition~\ref{prop:norms}. 
Relation \eqref{aen-1n3} implies that
\beq
E_{n,n-1}\xi_{\mu\nu}=\sum_{i=1}^{n-1} c_i(\mu,\nu)\ts \xi_{\mu-\delta_i,\nu}
\non\end{equation}
for some coefficients $c_i(\mu,\nu)$. Apply the operator $z_{j,n-1}$
to both sides of this relation. Since $z_{j,n-1}$ commutes with
$E_{n,n-1}$ we obtain from Lemma~\ref{lem:aximu} a recurrence relation
for the $c_i(\mu,\nu)$: if $\mu_j-\nu_j\geq 1$ then
\beq
c_i(\mu,\nu+\delta_j)=c_i(\mu,\nu)\ts \frac{m_i-\gamma_j-1}{m_i-\gamma_j},
\non\end{equation}
where $\gamma_j=\nu_j-j+1$. The proof is completed by induction.
The initial values of $c_i(\mu,\nu)$ are found by applying the relation
\beq
E_{n,n-1}\ts z_{n-1,i}=z_{ni}\frac{1}{h_i-h_{n-1}}+
z_{n-1,i}E_{n,n-1}\frac{h_i-h_{n-1}-1}{h_i-h_{n-1}}
\non\end{equation}
to the vector $\xi_{\mu}$ and taking into account that $E_{n,n-1}=z_{n,n-1}$.
Performing the calculation we get
\beq
E_{n,n-1}\ts \xi_{\mu\nu}=
\sum_{i=1}^{n-1} \frac{(m_i-\gamma_1)\cdots (m_i-\gamma_{n-2})}
{(m_i-m_1)\cdots\wedge_i\cdots (m_i-m_{n-1})}
\ts
\xi_{\mu-\delta_i,\nu}
\non\end{equation}
thus proving \eqref{aek+1k}.

\subsection{Characteristic identities}\label{subsec:ci}

Denote by $L$ the vector representation of $\gl_n$ and consider
its contragredient $L^*$. Note that $L^*$ is isomorphic to $L(0,\dots,0,-1)$.
Let $\{\ve_1,\dots,\ve_n\}$ denote the basis of $L^*$ dual to
the canonical basis $\{e_1,\dots,e_n\}$ of $L$.
Introduce the $n\times n$-matrix $E$ whose $ij$-th entry is the generator
$E_{ij}$. We shall interpret $E$ as the element
\beq
E=\sum_{i,j=1}^n e_{ij}\ot E_{ij}\in \End L^*\ot\U(\gl_n),
\non\end{equation}
where the $e_{ij}$ are the standard matrix units acting on $L^*$ by
$e_{ij}\ts \ve_k=\delta_{jk}\ts \ve_i$.	The basis element $E_{ij}$
of $\gl_n$ acts on $L^*$ as $-e_{ji}$ and hence $E$ may also be thought of as the image
of the element
\beq
e=-\sum_{i,j=1}^n E_{ji}\ot E_{ij}\in \U(\gl_n)\ot\U(\gl_n).
\non\end{equation}
On the other hand, using the standard coproduct $\Delta$ on $\U(\gl_n)$
defined by
\beq
\Delta(E_{ij})=E_{ij}\ot 1 + 1\ot E_{ij},
\non\end{equation}
we can write $e$ in the form
\beql{eimz}
e=\frac12\ts \Big(z\ot 1 +1\ot z -\Delta(z)\Big),
\end{equation}
where $z$ is the second order Casimir element
\beq
z=\sum_{i,j=1}^n E_{ij}E_{ji}\in \U(\gl_n).
\non\end{equation}

We have the tensor product decomposition
\beql{tensprll}
L^*\ot L(\la)\simeq L(\la-\delta_1)\oplus\cdots \oplus L(\la-\delta_n),
\end{equation}
where $L(\la-\delta_i)$ is considered to be zero if $\la_i=\la_{i+1}$.
On the level of characters this is a particular case of the {\it Pieri rule\/}
for the expansion of the product of a Schur polynomial by an 
elementary symmetric polynomial; see, e.g., Macdonald~\cite[Chapter~1]{m:sf}.
The Casimir element $z$
acts as a scalar operator in any highest weight representation $L(\la)$.
The corresponding eigenvalue is given by
\beq
z|^{}_{L(\la)}=\sum_{i=1}^n \la_i(\la_i+n-2i+1).
\non\end{equation}
Regarding now $E$ as an operator on $L^*\ot L(\la)$ and 
using \eqref{eimz} we derive
that the restriction of $E$ to the summand $L(\la-\delta_r)$ in \eqref{tensprll}
is the scalar operator with the eigenvalue $\la_r+n-r$ which we shall denote by $\al_r$.
This implies the {\it characteristic identity\/} for the matrix $E$,
\beql{charid}
\prod_{r=1}^n (E-\al_r)=0,
\end{equation}
as an operator in $L^*\ot L(\la)$. Moreover, the projection $P[r]$ of $L^*\ot L(\la)$
to the summand $L(\la-\delta_r)$ can be written explicitly as
\beq
P[r]=\frac{(E-\al_1)\cdots\wedge_r \cdots (E-\al_n)}
{(\al_r-\al_1)\cdots\wedge_r \cdots (\al_r-\al_n)}
\non\end{equation}
with $\wedge_r$ indicating that the $r$-th factor is omitted.
Together with \eqref{charid} this yields the {\it spectral decomposition\/} of $E$,
\beql{spectdec}
E=\sum_{r=1}^n\al_r\ts P[r].
\end{equation}
Consider the orthonormal Gelfand--Tsetlin bases
$\{\zeta^{}_{\La}\}$ of $L(\la)$ and
$\{\zeta^{}_{\La^{(r)}}\}$ of $L(\la-\delta_r)$ for $r=1,\dots,n$; see \eqref{orthonorm}.
Regarding the matrix element $P[r]_{ij}$ as an operator in $L(\la)$ we obtain
\beql{matelpr}
\lan \zeta^{}_{\La'}, P[r]_{ij}\ts \zeta^{}_{\La}\ran
= \lan \ve_i\ot\zeta^{}_{\La'}, P[r]\ts (\ve_j\ot\zeta^{}_{\La})\ran,
\end{equation}
where we have extended the inner products on $L^*$ and $L(\la)$ to $L^*\ot L(\la)$
by setting
\beq
\lan \eta\ot\zeta, \eta'\ot\zeta'\ran=\lan \eta, \eta'\ran
\lan \zeta, \zeta'\ran
\non\end{equation}
with $\eta,\eta'\in L^*$ and $\zeta, \zeta'\in L(\la)$.
Furthermore, using the expansions
\beq
\ve_j\ot\zeta^{}_{\La}=\sum_{s=1}^n\sum_{\La^{(s)}_{}} 
\lan \ve_j\ot\zeta^{}_{\La},\zeta^{}_{\La^{(s)}}\ran\ts\zeta^{}_{\La^{(s)}}, 
\non\end{equation}
bring \eqref{matelpr} to
the form
\beq
\sum_{\La^{(r)}_{}} 
\lan \ve_i\ot\zeta^{}_{\La'}, \zeta^{}_{\La^{(r)}}\ran
\lan \ve_j\ot\zeta^{}_{\La},  \zeta^{}_{\La^{(r)}}\ran,
\non\end{equation}
where we have used the fact that $P[r]$ is the identity map on
$L(\la-\delta_r)$, and zero on $L(\la-\delta_s)$ with $s\ne r$.
The numbers $\lan \ve_i\ot\zeta^{}_{\La'}, \zeta^{}_{\La^{(r)}}\ran$
are the {\it Wigner coefficients\/} (a particular case
of the {\it Clebsch--Gordan coefficients\/}). They can be used
to express the matrix elements of the generators $E_{ij}$ in the Gelfand--Tsetlin
basis as follows. Using the spectral decomposition \eqref{spectdec} we get
\beq
E_{ij}=\sum_{r=1}^n \al_r\ts P[r]_{ij}.
\non\end{equation}
Therefore, we derive the following result from \eqref{matelpr}. 

\bth\label{thm:wcme}
We have
\beq
\lan \zeta^{}_{\La'}, E_{ij}\ts \zeta^{}_{\La}\ran
=\sum_{r=1}^n \al_r\sum_{\La^{(r)}_{}} 
\lan \ve_i\ot\zeta^{}_{\La'}, \zeta^{}_{\La^{(r)}}\ran
\lan \ve_j\ot\zeta^{}_{\La},  \zeta^{}_{\La^{(r)}}\ran.
\non\end{equation}
\eth

Employing the characteristic identities for both the Lie algebras
$\gl_{n+1}$ and $\gl_n$ it is possible to determine the values
of the Wigner coefficients and thus to get an independent derivation
of the formulas of Theorem~\ref{thm:abasis}. In fact,
explicit formulas for the matrix elements of $E_{ij}$ with $|i-j|>1$ can also be given;
see Gould~\cite{g:me} for details.

The approach based upon the characteristic identities also leads
to an alternative presentation of the lowering and raising operators.
Taking $\zeta^{}_{\La}$ to be the highest vector $\xi$ in \eqref{matelpr}
we conclude that $P[r]_{ij}\ts \xi=0$ for $j>r$. Consider now $\gl_n$
as a subalgebra of $\gl_{n+1}$. Suppose that $\xi$ is a highest vector
of weight $\la$ in a representation $L(\la')$ of $\gl_{n+1}$. 
The previous observation implies that the vector
\beq
\sum_{i=r}^n E_{n+1,i}\ts P[r]_{ir}\ts\xi
\non\end{equation}
is again a $\gl_n$-highest vector of weight $\la-\delta_r$.

\bpr\label{prop:lozg}
We have the identity of operators 
on the space $L(\la')^+_{\la}$:
\beq
p\ts E_{n+1,r}=\sum_{i=r}^{n} E_{n+1,i}\ts P[r]_{ir}
\non\end{equation}
where $p$ is the extremal projector for $\gl_{n}$. 
\epr

\Outline  Since the both sides represent lowering operators they must be
proportional. It is therefore sufficient to apply both
sides to a vector $\xi\in L(\la')^+_{\la}$ and compare
the coefficients at $E_{n+1,r}\ts\xi$. For the calculation we use
the explicit formula \eqref{apni}  for  $p\ts E_{n+1,r}$
and the relation
\beq
P[r]_{rr}\ts\xi=\prod_{s=r+1}^n\frac{h_r-h_s-1}{h_r-h_s}\ts\xi
\non\end{equation}
which can be derived from the characteristic identities.
\endproof

An analogous argument leads to a similar formula for the raising operators.
Here one starts with the dual characteristic identity
\beq
\prod_{r=1}^n (\overline E-\bar\al_r)=0,
\non\end{equation}
where the $ij$-th matrix element of $\overline E$ is $-E_{ij}$,
$\bar\al_r=-\la_r+r-1$ and the powers of $\overline E$
are defined recursively by
\beq
(\overline E^{\ts p})_{ij}=\sum_{k=1}^n (\overline E^{\ts p-1})_{kj}\overline E_{ik}.
\non\end{equation}
For any $r=1,\dots,n$ the dual projection operator is given by
\beq
\overline P[r]=\frac{(\overline E-\bar\al_1)\cdots\wedge_r \cdots 
(\overline E-\bar\al_n)}
{(\bar\al_r-\bar\al_1)\cdots\wedge_r \cdots (\bar\al_r-\bar\al_n)}.
\non\end{equation}

\bpr\label{prop:rozg}
We have the identity of operators 
on the space $L(\la')^+_{\la}$:
\beq
p\ts E_{r,n+1}=\sum_{i=1}^{r}E_{i,n+1} \ts \overline P[r]_{ri}.
\non\end{equation}
\epr

\subsection{Quantum minors}\label{subsec:qm}

For a complex parameter $u$ introduce the $n\times n$-matrix $E(u)=u\ts 1 +E$.
Given sequences $a_1,\ldots, a_s$ and $b_1,\ldots, b_s$ of elements
of $\{1,\dots,n\}$ the corresponding {\it quantum minor\/} 
of the matrix $E(u)$
is defined
by the following equivalent formulas:
\begin{align}\label{qminor}
{E(u)}^{a_1\cdots\ts a_s}_{b_1\cdots\ts b_s}&=
\sum_{\sigma\in \Sym_s} \sgn \sigma \cdot E(u)_{a_{\sigma(1)}b_1}\cdots
E(u-s+1)_{a_{\sigma(s)}b_s}\\
\label{qminor2}
{}&=
\sum_{\sigma\in \Sym_s} \sgn \sigma\cdot E(u-s+1)_{a_1b_{\sigma(1)}}\cdots
E(u)_{a_sb_{\sigma(s)}}.
\end{align}
This is a polynomial in $u$ with coefficients in $\U(\gl_n)$. It
is skew symmetric under permutations of the indices $a_i$, or $b_i$.

For any index $1\leq i< n$ introduce the polynomials
\beq
\tau^{}_{ni}(u)=
{E(u)\ts}^{i+1\ts\cdots\ts n}_{i\ts\cdots\ts n-1}\qquad\text{and}\qquad
\tau^{}_{in}(u)=(-1)^{i-1}{E(u)\ts}^{1\cdots\ts i}_{1\cdots\ts i-1,n}.
\non\end{equation}
For instance,
\beq
\bal
\tau^{}_{13}(u)&=E_{13},\qquad \tau^{}_{23}(u)=-E_{23}\ts(u+E_{11})+E_{21}E_{13},\\
\tau^{}_{32}(u)&=E_{32},\qquad \tau^{}_{31}(u)=E_{21}E_{32}-E_{31}\ts(u+E_{22}-1).
\eal
\non\end{equation}

\bpr\label{prop:qmlr}
If $\eta\in L(\la)^+_{\mu}$ then
\beq
\tau^{}_{ni}(-\mu_i)\ts\eta\in L(\la)^+_{\mu-\delta_i}\qquad\text{and}\qquad
\tau^{}_{in}(-\mu_i+i-1)\ts\eta\in L(\la)^+_{\mu+\delta_i}.
\non\end{equation}
\epr

\Outline  The proof is based upon the
following relations
\beql{qmrel}
[E_{ij}, {E(u)\ts}^{a_1\cdots\ts a_s}_{b_1\cdots\ts b_s}]=
\sum_{r=1}^s\left(\delta_{ja_r}
{E(u)\ts}^{a_1\cdots\ts i\cdots\ts a_s}_{b_1\cdots\ts b_s}
-\delta_{ib_r}
{E(u)\ts}^{a_1\cdots\ts a_s}_{b_1\cdots\ts j\cdots\ts b_s}
\right),
\end{equation}
where $i$ and $j$ on the right hand side take the $r$-th positions.
\endproof

The relations \eqref{qmrel} imply the important
property of the quantum minors: for any indices $i,j$ we have
\beq
[E_{a_ib_j},
{E(u)\ts}^{a_1\cdots\ts a_s}_{b_1\cdots\ts b_s}]=0.
\non\end{equation}
In particular, this implies the centrality of the Capelli determinant
${\cal C}(u)={E(u)\ts}^{1\cdots\ts n}_{1\cdots\ts n}$; see \eqref{aqdet}.

The lowering and raising operators
of Proposition~\ref{prop:qmlr} can be shown to essentially coincide
with those defined in Section~\ref{subsec:loexp}. 
To write down the formulas we shall need to evaluate the variable $u$
in $\U(\h)$. To make this operation well-defined we accept the agreement
used in the evaluation of the Capelli determinant in \eqref{capeva}.

\bpr\label{prop:eqvoqm}
We have the identities for any $i=1,\dots,n-1$ 
\beql{taueqlz}
\tau^{}_{ni}(-h_i-i+1)=z_{ni}\qquad\text{and}\qquad
\tau^{}_{in}(-h_i)=z_{in}.
\end{equation}
\epr

Using this interpretation of the lowering operators
one can express the Gelfand--Tsetlin basis vector \eqref{axiL}
in terms of the quantum minors $\tau^{}_{ki}(u)$.
The action of certain other quantum minors
on these vectors can be explicitly found.
This will provide
one more independent proof of Theorem~\ref{thm:abasis}.
For $m\geq 1$ introduce the polynomials $A_m(u)$, $B_m(u)$ and $C_m(u)$ by
\beq
A_m(u)={E(u)\ts}^{1\cdots\ts m}_{1\cdots\ts m},\qquad
B_m(u)={E(u)\ts}^{1\cdots\ts m}_{1\cdots\ts m-1,m+1}, \qquad
C_m(u)={E(u)\ts}^{1\cdots\ts m-1,m+1}_{1\cdots\ts m}.
\non\end{equation}
We use the notation	 $l_{mi}=\la_{mi}-i+1$ and $l_i=\la_i-i+1$.

\bth\label{thm:dg} 
Let $\{\xiL\}$ be the Gelfand--Tsetlin basis of $L(\la)$. Then 
\begin{align}\label{agt}
A_m(u)\ts\xiL&=(u+l_{m1})\cdots (u+l_{mm})\ts\xiL,
\\
B_m(-l_{mj})\ts\xiL&=-\prod_{i=1}^{m+1} (l_{m+1,i}-l_{mj})
\ts \xi^{}_{\La+\delta_{mj}}\quad\text{for}\quad j=1,\dots,m, 
\non\\
\label{cgt}
C_m(-l_{mj})\ts\xiL&=\prod_{i=1}^{m-1} (l_{m-1,i}-l_{mj})
\ts \xi^{}_{\La-\delta_{mj}}\quad\text{for}\quad j=1,\dots,m,
\end{align}					   
where $\La\pm\delta_{mj}$ is obtained from
$\La$ by replacing the entry $\lambda_{mj}$ with $\lambda_{mj}\pm1$.
\endproof
\eth

Applying the Lagrange interpolation formula we can find 
the action
of $B_m(u)$ and $C_m(u)$ for any $u$. 
Note that these polynomials have degree $m-1$
with the leading coefficients $E_{m,m+1}$ and
$E_{m+1,m}$, respectively. 
Theorem~\ref{thm:abasis} is therefore an 
immediate corollary of Theorem~\ref{thm:dg}.

Formula \eqref{cgt} prompts a quite different construction of the basis vectors
of $L(\la)$ which uses the polynomials $C_m(u)$ instead of the
traditional lowering operators $z_{ni}$. Indeed, for a particular value of $u$,
$C_m(u)$ takes a basis vector into another one, up to a factor.
Given a pattern $\La$ associated with $\la$,
define the vector $\kappa^{}_{\La}\in L(\la)$ by
\begin{multline}\label{kappala}
\kappa^{}_{\La}=\prod_{k=1,\dots,n-1}^{\rightarrow}
\Big\{\ts C_{n-1}(-l_{n-1,k}-1)\cdots C_{n-1}(-l_k+1)\ts C_{n-1}(-l_k)\\
{}\times C_{n-2}(-l_{n-2,k}-1)\cdots C_{n-2}(-l_k+1)\ts C_{n-2}(-l_k)	\\
{}\times \cdots \times C_{k}(-l_{kk}-1)\cdots C_{k}(-l_k+1)\ts C_{k}(-l_k)\Big\}\ts\xi.
\non
\end{multline}

\bth\label{thm:ntgt}
The vectors $\kappa^{}_{\La}$ with $\La$ running over
all patterns associated with $\la$ form a basis
of $L(\la)$ and one has
$
\kappa^{}_{\La}=N^{}_{\La}\ts\xi^{}_{\La},
$
for a nonzero constant $N^{}_{\La}$.
\eth

The value of the constant $N^{}_{\La}$ can be found
from \eqref{cgt}. Using the relations between the elements
$A_m(u)$, $B_m(u)$ and $C_m(u)$ one can derive
Theorem~\ref{thm:dg} from Theorem~\ref{thm:ntgt}
with the use of Proposition~\ref{prop:eigenv} below; 
see Nazarov and Tarasov~\cite{nt:yg}
for details.

Observe that $A_m(u)$ is the Capelli determinant \eqref{aqdet} for the Lie algebra
$\gl_m$. Therefore, its coefficients $a_{mi}$ defined by
\beq
A_m(u)=u^m+a_{m1}\ts u^{m-1}+\cdots+a_{mm}
\non\end{equation}
are generators of the center of the enveloping algebra $\U(\gl_m)$.
All together the elements $a_{mi}$ with $1\leq i\leq m\leq n$
generate a commutative subalgebra ${\cal A}_n$ of $\U(\gl_n)$ which is called
the {\it Gelfand--Tsetlin subalgebra\/}. By \eqref{agt}, the basis vectors $\xi^{}_{\La}$
are simultaneous eigenvectors for the elements of the subalgebra ${\cal A}_n$.
Introduce the corresponding eigenvalues of the generators $a_{mi}$ by
\beql{ami}
a_{mi}\ts \xi^{}_{\La} =\al_{mi}(\La)\ts \xi^{}_{\La}.
\end{equation}
Thus, $\al_{mi}(\La)$ is the $i$-th 
elementary symmetric polynomial in $l_{m1},\dots,l_{mm}$.

\bpr\label{prop:eigenv}
For any pattern $\La$ associated with $\la$,
the one-dimensional subspace of $L(\la)$ spanned by the basis vector $\xi^{}_{\La}$
is uniquely determined by the set of eigenvalues $\{\al_{mi}(\La)\}$.
\epr

\subsection*{Bibliographical notes}\label{subsec:bibgln}


The explicit formulas for the lowering and raising operators \eqref{zin}
and \eqref{zni} first appeared in Nagel and Moshinsky~\cite{nm:ol};
see also Hou Pei-yu~\cite{h:ob} and Zhelobenko~\cite{z:cl}.
The derivation of the Gelfand--Tsetlin formulas outlined in Section~\ref{subsec:loexp}
follows
Zhelobenko~\cite{z:cl} and
Asherova, Smirnov and Tolstoy~\cite{ast:po2}. 
The extremal projectors were originally
introduced by Asherova, Smirnov and Tolstoy~\cite{ast:po} (see also \cite{ast:dc}).
In a subsequent paper \cite{ast:po2} the projectors were used to
construct the lowering operators and derive the relations
between them. A systematic study of the extremal projectors and the corresponding
Mickelsson algebras was undertaken by Zhelobenko: a detailed exposition
is given in his paper \cite{z:it} and book \cite{z:rr}. 
The application to the Gelfand--Tsetlin 
formulas is contained in his paper \cite{z:gz}.
Section~\ref{subsec:ma} is a brief outline of the general results
which are used in the basis constructions.

The first proof of Theorem~\ref{thm:mick} was given by 
van den Hombergh~\cite{h:nm}
as an answer to the question posed by Mickelsson~\cite{m:sa}.
A derivation of the relations in the Mickelsson--Zhelobenko 
algebra $\Z(\gl_n,\gl_m)$ with the use of the Capelli-type determinants
is contained in the author's paper \cite{m:yt}. 
A proof of the formulas \eqref{atzz} and
\eqref{atzz2} is also given there. 
The results of Section~\ref{subsec:ci}
are due to Gould~\cite{g:cir, g:ia, g:me, g:gr}.	
The characteristic identity \eqref{charid}
was proved by Green~\cite{gr:ci}. 
The significance of the Wigner coefficients in mathematical physics
is discussed in the book by Biedenharn and Louck~\cite{bl:am}.   
The definition \eqref{qminor}
of the quantum minors
is inspired by the theory of ``quantum" algebras
called the {\it Yangians\/}; see \cite{m:ya, mno:yc} for a review of the theory.
The polynomials $A_m(u)$, $B_m(u)$ and $C_m(u)$ are essentially the images
of the {\it Drinfeld generators\/} of the Yangian $\Y(n)$
under the evaluation homomorphism to the universal enveloping algebra
$\U(\gl_n)$. 
The quantum minor presentation of the lowering operators \eqref{taueqlz}
is due to the author \cite{m:gt}; see also \cite{m:yt}.
The construction of the Gelfand--Tsetlin basis vectors $\kappa^{}_{\La}$
with the use of the Drinfeld generators (Theorem~\ref{thm:ntgt})
was devised by Nazarov and Tarasov \cite{nt:yg}.

Analogs of the extremal projector were given by 
Tolstoy~\cite{t:ep, t:epc, t:epr, t:epq, t:po} 
for a wide class
of Lie (super)algebras and their quantized enveloping algebras.
The corresponding super and quantum versions of 
the Mickelsson--Zhelobenko algebras are studied in
\cite{t:epr, t:epq, t:po}. 
An alternative ``tensor formula" for the extremal projector
was provided by Tolstoy and Draayer~\cite{td:na}.
The techniques of extremal
projectors were applied by Khoroshkin and Tolstoy~\cite{kt:ep}
for calculation of the universal $R$-matrices for 
quantized enveloping algebras.
A basis of Gelfand--Tsetlin type for representations
of the exceptional Lie algebra $G_2$ was constructed
by Sviridov, Smirnov and Tolstoy~\cite{sst:cw, sst:si}.

Bases of Gelfand--Tsetlin type have been constructed for
representations of various types of algebras. For the quantized
enveloping algebra $\U_q(\gl_n)$ such bases were constructed
by Jimbo~\cite{j:qr}, Ueno, Takebayashi and  Shibukawa~\cite{uts:gz},
Nazarov and Tarasov~\cite{nt:yg}, Tolstoy~\cite{t:epq}.
The results of \cite{nt:yg}
include $q$-analogs of Theorems~\ref{thm:dg} and \ref{thm:ntgt},
while \cite{t:epq} contains matrix element formulas
for the generators corresponding to arbitrary roots.
Gould and Biedenharn~\cite{gb:pc} developed pattern calculus
for representations of the quantum group $\U_q(u(n))$.
Polynomial realizations of the Gelfand--Tsetlin basis
for representations of $\U_q(\sll_3)$ were given by
Dobrev and Truini~\cite{dt:ir, dt:pr} and Dobrev, Mitov and Truini~\cite{dmt:ng}.

Gelfand--Tsetlin bases
for `generic' representations of
the Yangian $\Y(n)$ were constructed in \cite{m:gt}.
Theorem~\ref{thm:dg} was proved there in a more general context
of representations of the {\it Yangian of level\/} $p$ for $\gl_n$ which was
previously introduced by Cherednik~\cite{c:qg}.	
In particular, the enveloping algebra $\U(\gl_n)$
coincides with the Yangian of level $1$.
A more general class of the {\it tame\/} 
Yangian modules was introduced and explicitly
constructed by Nazarov and Tarasov~\cite{nt:ry} via the {\it trapezium\/} or {\it skew\/}
analogs of the Gelfand--Tsetlin patterns.
Their approach was motivated by the so-called {\it centralizer construction\/}
devised by Olshanski~\cite{o:ea, o:yu, o:ri} and also 
employed by Cherednik~\cite{c:ni, c:qg}.
Basis vectors in the tame Yangian modules are characterized in a way
similar to Proposition~\ref{prop:eigenv}.
The skew Yangian modules were also studied in \cite{m:yt}
with the use of the quantum Sylvester theorem and the Mickelsson algebras. 

The center of $\U(\gl_n)$ possesses 
several natural families of generators and so does
the Gelfand--Tsetlin subalgebra ${\cal A}_n$.
The corresponding eigenvalues in $L(\la)$ are known explicitly;
see, e.g., \cite{m:ya} for a review.
An alternative description of ${\cal A}_n$
was given by Gelfand, Krob, Lascoux,
Leclerc,  Retakh and Thibon~\cite[Section~7.3]{gkllrt:ns} as an application
of their theory of noncommutative symmetric 
functions and quasi-determinants. 

The combinatorics of the skew Gelfand--Tsetlin patterns is employed
by Berenstein and Zelevinsky~\cite{bz:ig} to obtain 
multiplicity formulas for the skew representations of $\gl_n$.
Applications to continuous piecewise linear actions of the symmetric group
were found by Kirillov and Berenstein~\cite{kb:gg}.

The explicit realization of irreducible finite-dimensional
representations of $\gl_n$ via the Gelfand--Tsetlin bases
has important applications in the representation theory
of the quantum affine algebras and Yangians. In particular,
Theorem~\ref{thm:dg} and its Yangian analog \cite{m:gt}
are crucial in the proof of the irreducibility
criterion of the tensor products of the 
Yangian evaluation modules (a generalization to $\gl_n$ of
Theorem~\ref{thm:irrcr} below); see \cite{m:ic}.

Analogs of the Gelfand--Tsetlin bases for representations
of some Lie superalgebras and their
quantum analogs were given by Ottoson~\cite{o:cuo, o:cu},
Palev~\cite{p:fd1, p:fd2, p:if, p:et},
Palev, Stoilova and van der Jeugt~\cite{psj:fd},
Palev and Tolstoy~\cite{pt:fd},
Tolstoy, Istomina and Smirnov~\cite{tis:gt}.
Highest weight irreducible representations
for the Lie (super)algebras
of infinite matrices and their quantum analogs were
constructed by Palev~\cite{p:hw1, p:hw2}
and Palev and Stoilova~\cite{ps:hw1, ps:hw2, ps:hw3}
via bases of Gelfand--Tsetlin-type.  

The explicit formulas of Theorem~\ref{thm:abasis} make it possible
to define a class of infinite-dimensional representations of
$\gl_n$ by altering the inequalities \eqref{aconl}. Families of
such representations were introduced by Gelfand and Graev~\cite{gg:fd}.
However, as was later observed by Lemire and Patera~\cite{lp:fa},
some necessary conditions were missing in \cite{gg:fd} so that only a part
of those families actually provides representations. More general theory
of the so-called {\it Gelfand--Tsetlin modules\/} is developed 
by Drozd, Futorny and Ovsienko~\cite{dfo:iw, dfo:gz, dfo:gz3, dfo:hc},
Ovsienko~\cite{o:fs} and 
Mazorchuk~\cite{m:gv, m:cg}.
The starting point of the theory
is to axiomatize the property of the basis vectors \eqref{ami}
and to consider the module
generated by an eigenvector for the Gelfand--Tsetlin subalgebra
with a given arbitrary set of eigenvalues
$\{\al_{mi}\}$.
Some $q$-analogs of such modules were constructed by 
Mazorchuk and Turowska~\cite{mt:gz}.

The formulas of Theorem~\ref{thm:abasis} were applied by Olshanski~\cite{o:du, o:iu}
to study unitary representations of the pseudo-unitary
groups $U(p,q)$. In particular, he
classified all irreducible unitarizable highest weight representations
of the Lie algebra $\mathfrak{u}(p,q)$~\cite{o:du}. 
This work was extended by the author to a family of
the Enright--Varadarajan
modules over $\mathfrak{u}(p,q)$~\cite{m:us}.
Analogs of the Gelfand--Tsetlin bases for the unitary highest weight modules
were constructed in \cite{m:gtb}.

Applications of the Gelfand--Tsetlin bases 
in mathematical physics
are reviewed in the books by Barut and R\c aczka~\cite{br:tg}
and Biedenharn and Louck~\cite{bl:am}.


\section{Weight bases for representations of $\oa_N$ and $\spa_{2n}$}\label{sec:osp}
\setcounter{equation}{0}

Let $\g_n$ denote the rank $n$ simple complex Lie algebra of type
$B,C,$ or $D$. That is,
\beql{oos}
\g_n=\oa_{2n+1},\quad \spa_{2n},\quad\text{or}\quad\oa_{2n},
\end{equation}
respectively. Let	$V(\lambda)$ denote
the finite-dimensional
irreducible representation  of $\g_n$ 
with the highest weight $\lambda$. 
The restriction of $V(\lambda)$
to the subalgebra
$\g_{n-1}$ is not multiplicity-free in general.
This means that if $V'(\mu)$ is the finite-dimensional
irreducible	representation  of $\g_{n-1}$ 
with the highest weight $\mu$, then the space
\beql{homsp}
\Hom_{\g^{}_{n-1}}(V'(\mu),V(\la))
\end{equation}
need not be one-dimensional. In order to construct
a basis of $V(\lambda)$ associated with the chain
of subalgebras
\beq
\g_1\subset\g_2\subset\cdots\subset\g_n
\non\end{equation}
we need to construct a basis of the space \eqref{homsp}
which is isomorphic to the subspace $V(\lambda)^+_{\mu}$
of $\gl_{n-1}$-highest vectors of weight $\mu$ in $V(\la)$.
The subspace $V(\lambda)^+_{\mu}$
possesses a natural structure of a representation of the
centralizer $\C_n=\U(\g_n)^{\g_{n-1}}_{}$ of $\g_{n-1}$ in the
universal enveloping algebra $\U(\g_n)$.
It was shown by Olshanski~\cite{o:ty} that there exist
natural homomorphisms
\beq
\C_1
\leftarrow
\C_2
\leftarrow\cdots\leftarrow
\C_n\leftarrow
\C_{n+1}\leftarrow
\cdots.
\non\end{equation}
The projective limit of this chain turns out to be
an extension of the 
{\it twisted Yangian\/} $\Y^+(2)$ or $\Y^-(2)$, in the orthogonal
and symplectic case, respectively; see \cite{o:ty}, 
\cite{mno:yc} and \cite{mo:cc}
for the definition and properties of the twisted Yangians.
In particular, there is an algebra homomorphism
$\Y^{\pm}(2)\to \C_n$ which allows one to equip the space 
$V(\lambda)^+_{\mu}$	with a $\Y^{\pm}(2)$-module structure.
By the results of \cite{m:fd}, the representation
$V(\lambda)^+_{\mu}$
can be extended to a larger algebra, the {\it Yangian\/} $\Y(2)$.
This is a key fact which allows us to construct a natural basis in each
space $V(\lambda)^+_{\mu}$.
In the $C$ and $D$ cases the
$\Y(2)$-module $V(\lambda)^+_{\mu}$ is irreducible while
in the $B$ case it is a direct sum
of two irreducible submodules. This does not lead, however,
to major differences in the constructions, and the final
formulas are similar in all the three cases.

The calculations of the matrix elements of the generators of $\g_n$
are based on the relationship between
the twisted Yangian
$\Y^{\pm}(2)$ and the Mickelsson--Zhelobenko
algebra $\Z(\g_n,\g_{n-1})$; see Section~\ref{subsec:ma}.
We construct an algebra homomorphism $\Y^{\pm}(2)\to \Z(\g_n,\g_{n-1})$
which allows us to express the generators of the twisted Yangian,
as operators in $V(\lambda)^+_{\mu}$,
in terms of the lowering and raising operators.

\subsection{Raising and lowering operators}\label{subsec:rl}

Whenever possible we consider the three cases \eqref{oos} simultaneously,
unless otherwise stated.
The rows and columns of $2n\times 2n$-matrices will be enumerated
by the indices $-n,\dots,-1,1,\dots,n$, while 
the rows and columns of $(2n+1)\times (2n+1)$-matrices will be enumerated
by the indices $-n,\dots,-1,0,1,\dots,n$. Accordingly, the index $0$
will usually be skipped in the former case.
For $-n\leq i,j\leq n$ set
\beql{fij}
F_{ij}=E_{ij}-\theta_{ij}\ts E_{-j,-i}
\end{equation}
where the $E_{ij}$ are the standard matrix units, and
\beq\label{theta}
\theta_{ij}=\begin{cases}
1\qquad&\text{in the orthogonal case},\\
\sgn i\cdot\sgn j\qquad&\text{in the symplectic case}.
\end{cases}
\end{equation}
The matrices $F_{ij}$ span the
Lie algebra $\g_n$.
The subalgebra $\g_{n-1}$ is spanned by the elements \eqref{fij} with the
indices $i,j$ running over the set $\{-n+1,\dots,n-1\}$. 
Denote by $\h=\h_n$ the diagonal Cartan subalgebra in $\g_n$. 
The elements $F_{11},\dots,F_{nn}$ form a basis of $\h$. 

The finite-dimensional irreducible representations of $\g_n$
are in a one-to-one correspondence with $n$-tuples
$\lambda=(\lambda_1,\dots,\lambda_n)$ 
where the numbers $\lambda_i$ satisfy the conditions
\beq\label{lamco}
\lambda_i-\lambda_{i+1}\in \ZZ_+ \qquad\text{for}\quad i=1,\dots,n-1,
\end{equation}
and
\beq\label{lam1}
\bsp
-2\ts \lambda_1&\in\ZZ_+ \qquad\text{for}\quad \g_n=\oa_{2n+1},\\
-\lambda_1&\in\ZZ_+ \qquad\text{for}\quad \g_n=\spa_{2n},\\
-\lambda_1-\lambda_2&\in\ZZ_+ \qquad\text{for}\quad \g_n=\oa_{2n}.
\end{split}
\end{equation}
Such an $n$-tuple $\lambda$ is called the {\it highest weight\/}
of the corresponding representation which
we shall denote by $V(\lambda)$.
It contains a unique, up to a constant factor, nonzero vector $\xi$
(the {\it highest vector\/}) such that
\beq
\begin{aligned}
F_{ii}\ts\xi&=\lambda_i\ts\xi\qquad
&\text{for}&\quad i=1,\dots,n,\\
F_{ij}\ts\xi&=0\qquad
&\text{for}&\quad -n\leq i<j\leq n.
\end{aligned}
\non\end{equation}
Denote by $V(\lambda)^+$ the subspace of $\g_{n-1}$-highest vectors
in $V(\lambda)$:
\begin{equation}
V(\lambda)^+=\{\eta\in V(\lambda)\ |\ F_{ij}\ts \eta=0,
\qquad -n<i<j<n\}.
\non
\end{equation}
Given a $\g_{n-1}$-weight
$\mu=(\mu_1,\dots,\mu_{n-1})$ we denote by $V(\lambda)^+_{\mu}$
the corresponding weight subspace in $V(\lambda)^+$:
\begin{equation}
V(\lambda)^+_{\mu}=\{\eta\in V(\lambda)^+\ |\ F_{ii}\ts\eta=
\mu_i\ts\eta,\qquad i=1,\dots,n-1\}.
\non
\end{equation}

Consider the  Mickelsson--Zhelobenko algebra  $\Z(\g_n,\g_{n-1})$ 
introduced in Section~\ref{subsec:ma}.
Let $p=p_{n-1}$ be the extremal projector
for the Lie algebra $\g_{n-1}$.
It satisfies the conditions
\beq
F_{ij}\ts p=p\ts F_{ji}=0\qquad\text{for}\quad -n<i<j<n. 
\non\end{equation}
By Theorem~\ref{thm:pbwma}, the elements
\begin{equation}\label{gener}
F_{nn},\qquad pF_{ia},\qquad a=-n,n,\quad i=-n+1,\dots,n-1
\end{equation}
are generators of $\Z(\g_n,\g_{n-1})$ in the orthogonal case. In the symplectic case,
the algebra $\Z(\g_n,\g_{n-1})$ is generated by the elements \eqref{gener}
together with $F_{n,-n}$ and $F_{-n,n}$.
To write down explicit formulas for the generators,
introduce the numbers $\rho_i$, where $i=1,\dots,n$, by
\beq
\rho_i=\begin{cases}
-i+1/2\qquad&\text{for}\quad \g_n=\oa_{2n+1},\\
-i\qquad&\text{for}\quad \g_n=\spa_{2n},\\
-i+1\qquad&\text{for}\quad \g_n=\oa_{2n}.
\end{cases}
\non\end{equation}
The numbers $-\rho_i$ are coordinates of the half-sum of positive roots
with respect to the upper triangular Borel subalgebra.
Now set
\beq
f_i=F_{ii}+\rho_i,\qquad f_{-i}=-f_i
\non\end{equation}
for $i=1,\dots,n$. Moreover, in the case of $\oa_{2n+1}$ also set
$f_0=-1/2$. The generators $pF_{ia}$ can be given by a uniform expression
in all the three cases. Let $a\in\{-n,n\}$
and $i\in\{-n+1,\dots,n-1\}$. Then we have modulo the ideal $\J'$,
\beql{pFia}
pF_{ia}=F_{ia}+\sum_{i>i_1>\cdots>i_s>-n}
F_{ii_1}F_{i_1i_2}\cdots F_{i_{s-1}i_s}F_{i_sa}
\frac{1}{(f_i-f_{i_1})\cdots (f_i-f_{i_s})},
\end{equation}
summed over $s\geq 1$. It will be convenient to work with
normalized generators of $\Z(\g_n,\g_{n-1})$. 
Set
\beq
z_{ia}=pF_{ia}\ts (f_i-f_{i-1})\cdots(f_i-f_{-n+1})
\non\end{equation}
in the $B,C$ cases, and
\beq
z_{ia}=pF_{ia}(f_i-f_{i-1})\cdots\wh{(f_i-f_{-i})}\cdots(f_i-f_{-n+1})
\non\end{equation}
in the $D$ case, where
the hat indicates the factor to be omitted if it occurs.
We shall also use the elements $z_{ai}$ defined by
\beq
z_{ai}=(-1)^{n-i}\ts z_{-i,-a}\qquad\text{and}\qquad
z_{ai}=(-1)^{n-i}\ts\sgn a\cdot z_{-i,-a},
\non\end{equation}
in the orthogonal and
symplectic case, respectively. The elements $z_{ia}$ satisfy some quadratic
relations which can be shown to be the defining relations of the algebra
$\Z(\g_n,\g_{n-1})$. In particular, we have for all $a,b\in\{-n,n\}$
and $i+j\ne 0$,
\beql{bcdrelcom}
z_{ia}z_{jb}+z_{ja}z_{ib}(f_i-f_j-1)=z_{ib}z_{ja}(f_i-f_j).
\end{equation}
Thus, $z_{ia}$ and $z_{ja}$ commute for $i+j\ne 0$.
Also, $z_{ia}$ and $z_{ib}$ commute for
$i\ne 0$ and all values of $a$ and $b$.
Analogs of the relation \eqref{z_ins_ni} 
in the algebra $\Z(\g_n,\g_{n-1})$ can be explicitly written down as well.
However, we shall avoid using them in a way similar to the proof
of Lemma~\ref{lem:aximu}.

The elements $z_{ia}$ 
naturally act in the space 
$V(\lambda)^+$ by {\it raising\/} or {\it lowering\/} the weights.
We have for $i=1,\dots,n-1$:
\beq
z_{ia}:V(\lambda)^+_{\mu}\to V(\lambda)^+_{\mu+\delta_i},\qquad
z_{ai}:V(\lambda)^+_{\mu}\to V(\lambda)^+_{\mu-\delta_i},\non
\non\end{equation}
where $\mu\pm\delta_i$ is obtained from $\mu$ by replacing $\mu_i$
with $\mu_i\pm 1$. 
In the $B$ case the operators $z^{}_{0a}$ preserve each
subspace $V(\lambda)^+_{\mu}$. 

We shall need the following element which can be checked to belong
to the normalizer $\Norm \J'$, and so it can be regarded as
an element of the algebra $\Z(\g_n,\g_{n-1})$:
\beq
z_{n,-n}=
\sum_{n>i_1>\cdots>i_s>-n}
F_{ni_1}F_{i_1i_2}\cdots F_{i_s,-n}\ts
(f_{n}-f_{j_1})\cdots (f_{n}-f_{j_k})
\non\end{equation}
in the $B,C$ cases, and
\beq
z_{n,-n}=
\sum_{n>i_1>\cdots>i_s>-n}
F_{ni_1}F_{i_1i_2}\cdots F_{i_s,-n}\ts
\frac{(f_{n}-f_{j_1})\cdots (f_{n}-f_{j_k})}{2f_n}
\non\end{equation}
in the $D$ case,
where $s=0,1,\dots$ and $\{j_1,\dots,j_k\}$ is the complement
to the subset $\{i_1,\dots,i_s\}$ in
$\{-n+1,\dots,n-1\}$. 
The following is a counterpart of Lemma~\ref{lem:elow}
and is crucial in the calculation of the matrix elements
of the generators in the bases.

\ble\label{lem:bcdflow}
For $a\in\{-n,n\}$ we have
\beq
F_{n-1,a}=\sum_{i=-n+1}^{n-1}z_{n-1,i}\ts z_{ia}
\frac{1}{(f_i-f_{-n+1})\cdots\wedge_i\cdots(f_i-f_{n-1})}
\non\end{equation}
in the $B,C$ cases, and
\beq
F_{n-1,a}=\sum_{i=-n+1}^{n-1}z_{n-1,i}\ts z_{ia}
\frac{1}{(f_i-f_{-n+1})\cdots\wedge_{-i,i}\cdots(f_i-f_{n-1})}
\non\end{equation}
in the $D$ case,
where $z_{n-1,n-1}=1$ and the equalities are considered 
in $\U'(\g_n)$ modulo the ideal $\J'$. The wedge indicates
the indices to be skipped.
\ele

In order to write down the basis vectors, introduce the
interpolation polynomials $Z_{n,-n}(u)$ with coefficients
in the Mickelsson--Zhelobenko algebra $\Z(\g_n,\g_{n-1})$ by
\beql{bcinterp}
Z_{n,-n}(u)=\sum_{i=1}^{n}z_{ni}z_{i,-n}\prod_{j=1,\ts j\ne i}^{n}
\frac{u^2-g_j^2}{g_i^2-g_j^2}
\end{equation}
in the $B,C$ cases, and
\beql{dinterp}
Z_{n,-n}(u)=\sum_{i=1}^{n-1}z_{ni}z_{i,-n}\prod_{j=1,\ts j\ne i}^{n-1}
\frac{u^2-g_j^2}{g_i^2-g_j^2}
\end{equation}
in the $D$ case, where $g_i=f_i+1/2$.
Accordingly, we have
\beql{zeval}
Z_{n,-n}(g_i)=z_{ni}z_{i,-n}
\end{equation}
with the agreement that when $u$ is evaluated in $\U(\h)$, the coefficients
of the polynomial $Z_{n,-n}(u)$ are written to the left of the powers
of $u$, as appears in the formulas \eqref{bcinterp} and \eqref{dinterp}.

\subsection{Branching rules, patterns and basis vectors}\label{subsec:br}

The restriction of $V(\lambda)$ to the subalgebra $\g_{n-1}$
is given by
\beq
V(\lambda)|^{}_{\g_{n-1}}\simeq\underset{\mu}\oplus \ts c(\mu)\ts V'(\mu),
\non\end{equation}
where $V'(\mu)$ is the irreducible finite-dimensional representation of
$\g_{n-1}$ with the highest weight $\mu$.
The multiplicity $c(\mu)$ coincides with the dimension of
the space $V(\lambda)^+_{\mu}$, and its exact value is
found from the Zhelobenko branching rules~\cite{z:cg}.	We formulate them
separately for each case recalling the conditions \eqref{lamco} and \eqref{lam1}
on the highest weight~$\la$. In the formulas below we use the notation
\beq
l_i=\la_i+\rho_i+1/2,\qquad \gamma_i=\nu_i+\rho_i+1/2,
\non\end{equation}
where the $\nu_i$ are the parameters  defined in the branching rules.

A parameterization of basis vectors in $V(\lambda)$ is obtained by
applying the branching rules to its subsequent restrictions to 
the subalgebras of the chain
\beq
\g_1\subset\g_2\subset\cdots\subset\g_{n-1}\subset\g_n.
\non\end{equation}
This leads to the definition of the Gelfand--Tsetlin patterns for
the $B,C$ and $D$ types.
Then we give formulas for the basis vectors
of the representation $V(\lambda)$.
We use the notation
\beq
l_{ki}=\la_{ki}+\rho_i+1/2,\qquad l'_{ki}=\la'_{ki}+\rho_i+1/2,
\non\end{equation}
where the $\la_{ki}$ and $\la'_{ki}$ are 
the entries of the patterns defined below.

\bigskip
\noindent
{\bf B type case\/}. The multiplicity $c(\mu)$ 
equals the number
of $n$-tuples $(\nu^{\ts\prime}_1,\nu_2,\dots,\nu_n)$ satisfying the inequalities
\beq
\bal
&-\lambda_1\geq\nu^{\ts\prime}_1\geq\lambda_1\geq\nu_2\geq\lambda_2\geq \cdots\geq
\nu_{n-1}\geq\lambda_{n-1}\geq\nu_n\geq\lambda_n,\\
&-\mu_1\geq\nu^{\ts\prime}_1\geq\mu_1\geq\nu_2\geq\mu_2\geq \cdots\geq
\nu_{n-1}\geq\mu_{n-1}\geq\nu_n
\eal
\non\end{equation}
with $\nu^{\ts\prime}_1$ and all the $\nu_i$ being simultaneously
integers or half-integers together with the $\lambda_i$.
Equivalently, $c(\mu)$ 
equals the number of $(n+1)$-tuples $\nu=(\sigma,\nu_1,\dots,\nu_n)$,
with the entries given by
\beq
(\sigma,\nu_1)=\begin{cases}
(0,\nu^{\ts\prime}_1)\quad&\text{if $\nu^{\ts\prime}_1\leq 0$},\\
(1,-\nu^{\ts\prime}_1)\quad&\text{if $\nu^{\ts\prime}_1> 0$}.
\end{cases}
\non
\end{equation}

\ble\label{lem:bbasms}
The vectors
\beq
\xi_{\nu}=
z_{n0}^{\sigma}\ts\prod_{i=1}^{n-1}z_{ni}^{\nu_{i}-\mu_i}
z_{i,-n}^{\nu_{i}-\lambda_i}\cdot
\prod_{k=l_{n}}^{\gamma_{n}-1}Z_{n,-n}(k)\ts\xi
\non\end{equation}
form a basis of the space $V(\la)^+_{\mu}$.
\ele

Define the  {\it $B$ type 
pattern\/} $\Lambda$ associated with
$\lambda$ as an array of the form
\begin{align}
\sigma^{}_{n}\quad&\qquad\lambda^{}_{n1}\qquad\lambda^{}_{n2}
\qquad\qquad\cdots\qquad\qquad\lambda^{}_{nn}\non\\
&\lambda'_{n1}\qquad \lambda'_{n2}
\qquad\qquad\cdots\qquad\qquad\lambda'_{nn}\non\\
\sigma^{}_{n-1}\ &\qquad\lambda^{}_{n-1,1}\qquad\cdots
\qquad\lambda^{}_{n-1,n-1}\non\\
&\lambda'_{n-1,1}
\qquad\cdots\qquad\lambda'_{n-1,n-1}\non\\
\qquad\cdots&\qquad\cdots\qquad\cdots\non\\
\sigma^{}_{1}\quad&\qquad\lambda^{}_{11}\non\\
&\lambda'_{11}\non
\end{align}
such that $\lambda=(\lambda^{}_{n1},\dots, \lambda^{}_{nn})$,
each $\sigma^{}_{k}$ is $0$ or $1$, the remaining
entries are all non-positive
integers or non-positive half-integers together with the
$\lambda_i$, and the following inequalities hold
\beq
\lambda'_{k1}\geq\lambda^{}_{k1}\geq\lambda'_{k2}\geq
\lambda^{}_{k2}\geq \cdots\geq
\lambda'_{k,k-1}\geq\lambda^{}_{k,k-1}\geq\lambda'_{kk}\geq\lambda^{}_{kk}
\non
\end{equation}
for $k=1,\dots,n$, and
\beq
\lambda'_{k1}\geq\lambda^{}_{k-1,1}\geq\lambda'_{k2}\geq
\lambda^{}_{k-1,2}\geq \cdots\geq
\lambda'_{k,k-1}\geq\lambda^{}_{k-1,k-1}\geq\lambda'_{kk}
\non
\end{equation}
for $k=2,\dots,n$. In addition, in the case of the integer $\lambda_i$
the condition
\beq
\lambda'_{k1}\leq -1\qquad\text{if}\quad \sigma^{}_{k}=1
\non
\end{equation}
should hold for all $k=1,\dots,n$.

\bth\label{thm:bbasis}
The vectors
\beq
\xi^{}_{\Lambda}=
\prod_{k=1,\dots,n}^{\rightarrow}
\left(z_{k0}^{\sigma^{}_{k}}\cdot\prod_{i=1}^{k-1}
z_{ki}^{\lambda'_{ki}-\lambda^{}_{k-1,i}}
z_{i,-k}^{\lambda'_{ki}-\lambda^{}_{ki}}
\cdot\prod_{j=l^{}_{kk}}^{l'_{kk}-1}
Z_{k,-k}(j)\right)\xi
\non\end{equation}
parametrized by the 
patterns $\Lambda$ form a basis of 
the representation $V(\lambda)$. 
\eth

\bigskip
\noindent
{\bf C type case\/}. The multiplicity $c(\mu)$ 
equals the number
of $n$-tuples of integers $\nu=(\nu_1,\dots,\nu_n)$ satisfying the
inequalities
\beql{cineq}
\bal
&0\geq\nu_1\geq\lambda_1\geq\nu_2\geq\lambda_2\geq \cdots\geq
\nu_{n-1}\geq\lambda_{n-1}\geq\nu_n\geq\lambda_n,\\
&0\geq\nu_1\geq\mu_1\geq\nu_2\geq\mu_2\geq \cdots\geq
\nu_{n-1}\geq\mu_{n-1}\geq\nu_n.
\eal
\end{equation}

\ble\label{lem:cbasms}
The vectors
\beq
\xi_{\nu}=
\prod_{i=1}^{n-1}z_{ni}^{\nu_{i}-\mu_i}
z_{i,-n}^{\nu_{i}-\lambda_i}\cdot
\prod_{k=l_{n}}^{\gamma_{n}-1}Z_{n,-n}(k)\ts\xi
\non\end{equation}
form a basis of the space $V(\la)^+_{\mu}$.
\ele

Define the  {\it $C$ type 
pattern\/} $\Lambda$ associated with
$\lambda$ as an array of the form
\begin{align}
\quad&\qquad\lambda^{}_{n1}\qquad\lambda^{}_{n2}
\qquad\qquad\cdots\qquad\qquad\lambda^{}_{nn}\non\\
&\lambda'_{n1}\qquad \lambda'_{n2}
\qquad\qquad\cdots\qquad\qquad\lambda'_{nn}\non\\
\ &\qquad\lambda^{}_{n-1,1}\qquad\cdots
\qquad\lambda^{}_{n-1,n-1}\non\\
&\lambda'_{n-1,1}
\qquad\cdots\qquad\lambda'_{n-1,n-1}\non\\
&\qquad\cdots\qquad\cdots\non\\
\quad&\qquad\lambda^{}_{11}\non\\
&\lambda'_{11}\non
\end{align}
such that $\lambda=(\lambda^{}_{n1},\dots, \lambda^{}_{nn})$,
the remaining
entries are all non-positive
integers and the following inequalities hold
\beq
0\geq\lambda'_{k1}\geq\lambda^{}_{k1}\geq\lambda'_{k2}\geq
\lambda^{}_{k2}\geq \cdots\geq
\lambda'_{k,k-1}\geq\lambda^{}_{k,k-1}\geq\lambda'_{kk}\geq\lambda^{}_{kk}
\non
\end{equation}
for $k=1,\dots,n$, and
\beq
0\geq\lambda'_{k1}\geq\lambda^{}_{k-1,1}\geq\lambda'_{k2}\geq
\lambda^{}_{k-1,2}\geq \cdots\geq
\lambda'_{k,k-1}\geq\lambda^{}_{k-1,k-1}\geq\lambda'_{kk}
\non
\end{equation}
for $k=2,\dots,n$.

\bth\label{thm:cbasis}
The vectors
\beq
\xi^{}_{\Lambda}=
\prod_{k=1,\dots,n}^{\rightarrow}
\left(\prod_{i=1}^{k-1}
z_{ki}^{\lambda'_{ki}-\lambda^{}_{k-1,i}}
z_{i,-k}^{\lambda'_{ki}-\lambda^{}_{ki}}
\cdot\prod_{j=l^{}_{kk}}^{l'_{kk}-1}
Z_{k,-k}(j)\right)\xi
\non\end{equation}
parametrized by the 
patterns $\Lambda$ form a basis of 
the representation $V(\lambda)$. 
\eth

\bigskip
\noindent
{\bf D type case\/}. The multiplicity $c(\mu)$ 
equals the number of $(n-1)$-tuples $\nu=(\nu_1,\dots,\nu_{n-1})$
satisfying the inequalities
\beq
\bal
-|\lambda_1|&\geq\nu_1\geq\lambda_2\geq\nu_2\geq\lambda_3\geq
\cdots \geq\lambda_{n-1}\geq\nu_{n-1}\geq\lambda_n,\\
-|\mu_1|&\geq\nu_1\geq\mu_2\geq\nu_2\geq\mu_3\geq
\cdots \geq\mu_{n-1}\geq\nu_{n-1}
\eal
\non\end{equation}
with all the $\nu_i$ being simultaneously
integers or half-integers together with the $\lambda_i$.
Set $\nu_0=\text{max}\{\la_1,\mu_1\}$.

\ble\label{lem:dbasms}
The vectors
\beq
\xi_{\nu}=
\prod_{i=1}^{n-1}z_{ni}^{\nu_{i-1}-\mu_i}
z_{i,-n}^{\nu_{i-1}-\lambda_i}\cdot
\prod_{k=l_{n}}^{\gamma_{n-1}-2}Z_{n,-n}(k)\ts\xi
\non\end{equation}
form a basis of the space $V(\la)^+_{\mu}$.
\ele

Define the  {\it $D$ type 
pattern\/} $\Lambda$ associated with
$\lambda$ as an array of the form
\begin{align}
&\lambda^{}_{n1}\qquad\lambda^{}_{n2}
\qquad\qquad\qquad\cdots\qquad\qquad\qquad\lambda^{}_{nn}\non\\
&\qquad \lambda'_{n-1,1}
\qquad\qquad\cdots\qquad\qquad\lambda'_{n-1,n-1}\non\\
&\lambda^{}_{n-1,1}\qquad\cdots\qquad\lambda^{}_{n-1,n-1}\non\\
&\qquad\qquad\cdots\qquad\cdots\non\\
&\lambda^{}_{21}\qquad\lambda^{}_{22}\non\\
&\qquad\lambda'_{11}\non\\
&\lambda^{}_{11}\non
\end{align}
such that $\lambda=(\lambda^{}_{n1},\dots, \lambda^{}_{nn})$,
the remaining
entries are all non-positive
integers or non-positive half-integers together with the
$\lambda_i$, and the following inequalities hold
\begin{align}
-|\lambda^{}_{k1}|&\geq\lambda'_{k-1,1}\geq\lambda^{}_{k2}\geq
\lambda'_{k-1,2}\geq \cdots\geq
\lambda^{}_{k,k-1}\geq\lambda'_{k-1,k-1}\geq\lambda^{}_{kk},\non\\
-|\lambda^{}_{k-1,1}|&\geq\lambda'_{k-1,1}\geq\lambda^{}_{k-1,2}\geq
\lambda'_{k-1,2}\geq \cdots\geq
\lambda^{}_{k-1,k-1}\geq\lambda'_{k-1,k-1}\non
\end{align}
for $k=2,\dots,n$. Set $\la'_{k-1,0}=\text{max}\{\la^{}_{k1},\la^{}_{k-1,1}\}$.

\bth\label{thm:dbasis}
The vectors
\beq
\xi^{}_{\Lambda}=
\prod_{k=2,\dots,n}^{\rightarrow}
\left(\prod_{i=1}^{k-1}
z_{ki}^{\lambda'_{k-1,i-1}-\lambda^{}_{k-1,i}}
z_{i,-k}^{\lambda'_{k-1,i-1}-\lambda^{}_{ki}}
\cdot\prod_{j=l^{}_{kk}}^{l'_{k-1,k-1}-2}
Z_{k,-k}(j)\right)\xi
\non\end{equation}
parametrized by the 
patterns $\Lambda$ form a basis of 
the representation $V(\lambda)$. 
\eth

Proofs of Theorems~\ref{thm:bbasis}, \ref{thm:cbasis}
and \ref{thm:dbasis} will be outlined in the
next two sections. These are based on the application of 
the representation theory of the
twisted Yangians. Clearly, due to the branching rules,
it is sufficient to construct a basis in the multiplicity space $V(\lambda)^+_{\mu}$.

\subsection{Yangians and their representations}\label{subsec:yr}

We start by introducing the {\it Yangian\/}
$\Y(2)$ for the Lie algebra $\gl_2$. 
In what follows it will be convenient to 
use the indices $-n,n$ to
enumerate the rows 
and columns of $2\times 2$-matrices.
The Yangian $\Y(2)$ is the complex associative algebra with the
generators $t_{ab}^{(1)},t_{ab}^{(2)},\dots$ where 
$a,b\in\{-n,n\}$,
and the defining relations
\begin{equation}\label{rel}
(u-v)\ts [t_{ab}(u),t_{cd}(v)]=
t_{cb}(u)t_{ad}(v)-t_{cb}(v)t_{ad}(u),
\end{equation}
where
\beq
t_{ab}(u) = \delta_{ab} + t^{(1)}_{ab} u^{-1} + t^{(2)}_{ab}u^{-2} +
\cdots \in \Y(2)[[u^{-1}]].
\non\end{equation}
Introduce the series $s_{ab}(u)$, $a,b\in\{-n,n\}$ by
\beql{sab}
s_{ab}(u)=\theta_{nb}\ts t_{an}(u)t_{-b,-n}(-u)+
\theta_{-n,b}\ts t_{a,-n}(u)t_{-b,n}(-u)
\end{equation}
with $\theta_{ij}$ defined in \eqref{theta}.
Write
\beq
s_{ab}(u)=\delta_{ab}+s_{ab}^{(1)}u^{-1}+s_{ab}^{(2)}u^{-2}+\cdots.
\non\end{equation}
The {\it twisted Yangian\/} $\Y^{\pm}(2)$ is defined as the subalgebra of $\Y(2)$
generated by the elements $s_{ab}^{(1)},s_{ab}^{(2)},\dots$ where 
$a,b\in\{-n,n\}$. Also, $\Y^{\pm}(2)$ can be viewed as
an abstract algebra with generators $s_{ab}^{(r)}$ and quadratic
and linear defining relations
which have the following form
\beq
\bal
(u^2-v^2)\ts
[s_{ab}(u),s_{cd}(v)]=(u+v)\ts&\big(s_{cb}(u)s_{ad}(v)-s_{cb}(v)s_{ad}(u)\big)\\
{}-(u-v)\ts&\big(\theta_{c,-b}s_{a,-c}(u)s_{-b,d}(v)-
\theta_{a,-d}s_{c,-a}(v)s_{-d,b}(u)\big)\\
{}+\theta_{a,-b}&\big(s_{c,-a}(u)s_{-b,d}(v)-s_{c,-a}(v)s_{-b,d}(u)\big)
\eal
\non\end{equation}
and
\beq
\theta_{ab}\ts s_{-b,-a}(-u)=s_{ab}(u)\pm \frac{s_{ab}(u)-s_{ab}(-u)}{2u}. 
\non\end{equation}
Whenever the double sign $\pm{}$ or $\mp{}$ occurs,
the upper sign corresponds to the orthogonal case and the lower sign to
the symplectic case. In particular, we have the relation
\beq
[s_{n,-n}(u), s_{n,-n}(v)]=0.
\non\end{equation}

The Yangian $\Y(2)$ is a Hopf algebra with the 
coproduct 
\begin{equation}\label{cop}
\Delta (t_{ab}(u))=t_{an}(u)\ot
t_{nb}(u)+t_{a,-n}(u)\ot t_{-n,b}(u).
\end{equation}
The twisted Yangian $\Y^{\pm}(2)$ is a left coideal in $\Y(2)$
with
\begin{equation}\label{cops}
\Delta (s_{ab}(u))=\sum_{c,d\in\{-n,n\}}\theta_{bd}\ts 
t_{ac}(u)t_{-b,-d}(-u)\ot
s_{cd}(u).
\end{equation}

Given a pair of complex numbers $(\alpha,\beta)$ 
such that $\alpha-\beta\in\ZZ_+$
we denote by 
$L(\alpha,\beta)$ the irreducible representation of the Lie algebra
$\gl_2$ with the highest weight $(\alpha,\beta)$ with respect to the
upper triangular Borel subalgebra. 
Then $\dim L(\alpha,\beta)=\alpha-\beta+1$. We equip $L(\alpha,\beta)$ with
a $\Y(2)$-module structure by using the algebra homomorphism $\Y(2)\to\U(\gl_2)$
given by
\begin{equation}
t_{ab}(u)\mapsto \delta_{ab}+E_{ab}u^{-1},\qquad a,b\in\{-n,n\}.
\non
\end{equation}
The coproduct \eqref{cop} allows us to construct representations of
$\Y(2)$ of the form
\beql{tenpr}
L=L(\alpha_1,\beta_1)\ot\cdots\ot L(\alpha_k,\beta_k).
\end{equation}
Any finite-dimensional irreducible $\Y(2)$-module is isomorphic
to a representation of this type twisted by an automorphism
of $\Y(2)$ of the form
\beq
t_{ab}(u)\mapsto (1+\varphi_1\ts u^{-1}+\varphi_2\ts u^{-2}+\cdots)
\ts t_{ab}(u),\qquad \varphi_i\in\CC.
\non\end{equation}
There is an explicit irreducibility criterion for the $\Y(2)$-module $L$.
To formulate
the result, with each
$L(\alpha,\beta)$ associate the {\it string\/}
\beq
S(\alpha,\beta)=\{\beta,\beta+1,\dots,\alpha-1\}\subset\CC.
\non\end{equation}
We say that two strings $S_1$ and $S_2$ are 
{\it in general position\/}
if
\beq
\text{either}\quad
S_1\cup S_2\quad\text{is not a string, or}\quad
S_1\subseteq S_2,\quad\text{or}\quad S_2\subseteq S_1.
\non\end{equation}

\bth\label{thm:irrcr}
The representation \eqref{tenpr} of $\Y(2)$ is irreducible 
if and only if the strings
$S(\alpha_i,\beta_i)$, $i=1,\dots,k$, are pairwise in
general position.
\eth

Note that the generators $t_{ab}^{(r)}$ with $r>k$ act as zero operators
in $L$. Therefore, the operators $T_{ab}(u)=u^k\ts t_{ab}(u)$
are polynomials in $u$:
\begin{equation}\label{Tab}
T_{ab}(u)=\delta_{ab}\ts u^k+t_{ab}^{(1)}u^{k-1}+\cdots+t_{ab}^{(k)}.
\end{equation}
Let $\xi_i$ denote the highest vector of the
$\gl_2$-module $L(\alpha_i,\beta_i)$. 
Suppose that the $\Y(2)$-module $L$ given by \eqref{tenpr}
is irreducible and the strings $S(\alpha_i,\beta_i)$
are pairwise disjoint.
Set
\beql{tensprxi}
\eta=\xi_1\ot\cdots\ot\xi_k. 
\end{equation}
Then using \eqref{cop}
we easily check that $\eta$ is the highest vector
of the $\Y(2)$-module $L$. That is,	$\eta$ 
is annihilated by $T_{-n,n}(u)$, and it
is an eigenvector for the operators $T_{nn}(u)$ and 
$T_{-n,-n}(u)$. Explicitly,
\beq\label{eigen}
\bal
T_{-n,-n}(u)\ts\eta&=(u+\alpha_1)\cdots (u+\alpha_k)\ts\eta,\\
T_{nn}(u)\ts\eta&=(u+\beta_1)\cdots (u+\beta_k)\ts\eta.
\eal
\end{equation}
Let a $k$-tuple $\gamma=(\gamma_1,\dots,\gamma_k)$ satisfy
the conditions: for each $i$
\beq\label{coga}
\alpha_i-\gamma_i\in\ZZ_+,\qquad \gamma_i-\beta_i\in\ZZ_+.
\end{equation}
Set
\beq
\eta_{\gamma}=
\prod_{i=1}^{k}T_{n,-n}(-\gamma_i+1)\cdots 
T_{n,-n}(-\beta_i-1)T_{n,-n}(-\beta_i)\ts 
\eta.
\non\end{equation}
The following theorem provides a Gelfand--Tsetlin type
basis for representations of the Yangian $\Y(2)$ associated with the
embedding $\Y(1)\subset\Y(2)$. Here $\Y(1)$ is the 
(commutative) subalgebra of $\Y(2)$
generated by the elements $t_{nn}^{(r)}$, $r\geq 1$.

\bth\label{thm:basy}
Let the $\Y(2)$-module $L$ given by \eqref{tenpr}
be irreducible and the strings $S(\alpha_i,\beta_i)$
be pairwise disjoint. Then
the vectors $\eta_{\gamma}$ with $\gamma$ satisfying \eqref{coga}
form a basis of $L$. Moreover, the generators of $\Y(2)$ act
in this basis by the rule
\beq\label{actt}
\begin{aligned}
T_{nn}(u)\ts{\eta}_{\gamma}&=(u+\gamma_1)\cdots (u+\gamma_k)
\ts{\eta}_{\gamma}, \\
T_{n,-n}(-\gamma_i)\ts 
{\eta}_{\gamma}&={\eta}_{\gamma+\delta_i},\\
T_{-n,n}(-\gamma_i)\ts {\eta}_{\gamma}&=
-\prod_{m=1}^k(\alpha_m-\gamma_i+1)(\beta_m-\gamma_i)
\ts{\eta}_{\gamma-\delta_i},\\
T_{-n,-n}(u)\ts \eta_{\gamma}&=\prod_{i=1}^{k}
\frac{(u+\alpha_i+1)(u+\beta_i)}{u+\gamma_i+1}\ts \eta_{\gamma}\\
{}&+\prod_{i=1}^{k}
\frac{1}{u+\gamma_i+1}\ts T_{-n,n}(u)T_{n,-n}(u+1)\ts \eta_{\gamma}.
\eal
\end{equation}
\eth

These formulas are derived from the defining relations
for the Yangian \eqref{rel} with the use
of the {\it quantum determinant\/} 
\begin{align}\label{qdet1}
d(u)&=T_{-n,-n}(u+1)\ts T_{nn}(u)-T_{n,-n}(u+1)\ts T_{-n,n}(u)\\
\label{qdet2}
{}&=T_{-n,-n}(u)\ts T_{nn}(u+1)-T_{-n,n}(u)\ts T_{n,-n}(u+1).
\end{align}
The coefficients of the quantum determinant belong to the center of $\Y(2)$
and so, $d(u)$ acts in $L$ as a scalar which
can be found by the application of \eqref{qdet1} to the 
highest vector $\eta$. Indeed, 
by \eqref{eigen}
\begin{equation}
d(u)\ts \eta=(u+\alpha_1+1)\cdots (u+\alpha_{k}+1)
(u+\beta_1)\cdots(u+\beta_{k})\ts \eta.
\non
\end{equation}
This allows us to derive the last formula in \eqref{actt}
from \eqref{qdet2}.
The operators $T_{-n,n}(u)$ and $T_{n,-n}(u)$ are polynomials in $u$
of degree ${}\leq k-1$; see \eqref{Tab}. 
Therefore, their action can be found from
\eqref{actt} by using the Lagrange interpolation formula.

We can regard \eqref{tenpr}	as a module over the twisted Yangian $\Y^-(2)$
obtained by restriction. Irreducibility criterion of such a module
is provided by the following theorem.

\bth\label{thm:symirr}
The representation \eqref{tenpr} of $\Y^-(2)$ is irreducible 
if and only if each pair of strings
\beq%
S(\alpha_i,\beta_i),\  S(\alpha_j,\beta_j)
\qquad\text{and}\qquad  S(\alpha_i,\beta_i),\  S(-\beta_j,-\alpha_j)
\non\end{equation}
is in general position for all $i<j$.
\eth

The defining relations \eqref{rel} allow us
to rewrite the formula \eqref{sab} for $s_{n,-n}(u)$ in the form
\beq
s_{n,-n}(u)=\frac{u+1/2}{u}\Big(t_{n,-n}(u)t_{nn}(-u)-
t_{n,-n}(-u)t_{nn}(u)\Big).
\non\end{equation}
Therefore the operator in $L$ defined by
\begin{multline}\label{STT}
S_{n,-n}(u)=\frac{u^{2k}}{u+1/2}s_{n,-n}(u)=\\
\frac{(-1)^k}{u}\Big(T_{n,-n}(u)T_{nn}(-u)-
T_{n,-n}(-u)T_{nn}(u)\Big)
\end{multline}
is an even polynomial in $u$ of degree $\leq 2k-2$. 
Its action in the basis of $L$
provided in Theorem~\ref{thm:basy} is given by
\beq
S_{n,-n}(\gamma_i)\ts{\eta}_{\ts\gamma}=2\ts 
\prod_{a=1,\ts a\ne i}^k (-\gamma_i-\gamma_a) \ts 
{\eta}_{\ts\gamma+\delta_i},\qquad i=1,\dots,k.
\non\end{equation}
We have thus proved the following corollary.

\bco\label{cor:bassgen}
Suppose that the $\Y^-(2)$-module $L$ is irreducible and we have
\beq%
S(\alpha_i,\beta_i)\cap S(\alpha_j,\beta_j)=\emptyset
\qquad\text{and}\qquad  S(\alpha_i,\beta_i)\cap S(-\beta_j,-\alpha_j)=\emptyset
\non\end{equation}
for all $i<j$.\footnote{The second condition was erroneously omitted
in the formulation of \cite[Proposition~4.2]{m:br} although 
it is implicit in the proof.}
Then the vectors
\beq
\xi_{\gamma}=\prod_{i=1}^k S_{n,-n}(\gamma_i-1)\cdots 
S_{n,-n}(\beta_i+1)\ts S_{n,-n}(\beta_i)\ts\eta
\non\end{equation}
with $\gamma$ satisfying \eqref{coga} form a basis of $L$.
\eco

Let us now turn to the orthogonal twisted Yangian $\Y^+(2)$.
For any $\delta\in\CC$ denote by $W(\delta)$ the one-dimensional
representation of
$\Y^+(2)$ spanned by a vector $w$ such that
\begin{equation}
s_{nn}(u)\ts w=\frac{u+\delta}{u+1/2}\ts w,\qquad 
s_{-n,-n}(u)\ts w=\frac{u-\delta+1}{u+1/2}\ts w,
\non
\end{equation}
and $s_{a,-a}(u)\ts w=0$ for $a=-n,n$. 
By \eqref{cops} we can regard the tensor
product $L\ot W(\delta)$ as a representation of $\Y^+(2)$.
The representations of $\Y^+(2)$ of this type,
and the representations of $\Y^-(2)$ of type 
\eqref{tenpr}
essentially exhaust all finite-dimensional irreducible representations of
$\Y^{\pm}(2)$ \cite{m:fd}.

The following is an analog of Theorem~\ref{thm:symirr}.
\bth\label{thm:orthirr}
The representation $L\ot W(\delta)$ of $\Y^+(2)$ is irreducible 
if and only if each pair of strings
\beq%
S(\alpha_i,\beta_i),\  S(\alpha_j,\beta_j)
\qquad\text{and}\qquad  S(\alpha_i,\beta_i),\  S(-\beta_j,-\alpha_j)
\non\end{equation}
is in general position for all $i<j$, and  
none of the strings $S(\alpha_i,\beta_i)$ or $S(-\beta_i,-\alpha_i)$
contains $-\delta$.
\eth

Using the vector space isomorphism
\beql{isomL}
L\ot W(\delta)\to L,\qquad v\ot w\mapsto v,\qquad v\in L
\end{equation}
we can regard $L$ as a $\Y^+(2)$-module.
Accordingly, using
the defining relations \eqref{rel} and
the coproduct formula \eqref{cops} we can write
$s_{n,-n}(u)$, as an operator in $L$, in the form
\beq
s_{n,-n}(u)=\frac{u-\delta}{u}\ts t_{n,-n}(u)t_{nn}(-u)+
\frac{u+\delta}{u}\ts t_{n,-n}(-u)t_{nn}(u).
\non\end{equation}
Therefore the operator in $L$ defined by
\begin{multline}\label{orSTT}
S_{n,-n}(u)=u^{2k}\ts s_{n,-n}(u)=\\
\frac{(-1)^k}{u}\Big((u-\delta)\ts T_{n,-n}(u)T_{nn}(-u)+
(u+\delta)\ts T_{n,-n}(-u)T_{nn}(u)\Big)
\end{multline}
is an even polynomial in $u$ of degree $\leq 2k-2$. Its action in the basis of $L$
provided in Theorem~\ref{thm:basy} is given by
\beq
S_{n,-n}(\gamma_i)\ts{\eta}_{\ts\gamma}=2(-\delta-\gamma_i)\ts 
\prod_{a=1,\ts a\ne i}^k (-\gamma_i-\gamma_a) \ts 
{\eta}_{\ts\gamma+\delta_i},\qquad i=1,\dots,k.
\non\end{equation}
We have thus proved the following corollary.

\bco\label{cor:orbassgen}
Suppose that the $\Y^+(2)$-module $L\ot W(\delta)$ is irreducible
and we have
\beq%
S(\alpha_i,\beta_i)\cap S(\alpha_j,\beta_j)=\emptyset
\qquad\text{and}\qquad  S(\alpha_i,\beta_i)\cap S(-\beta_j,-\alpha_j)=\emptyset
\non\end{equation}
for all $i<j$. Then the vectors
\beq
\xi_{\gamma}=\prod_{i=1}^k S_{n,-n}(\gamma_i-1)\cdots 
S_{n,-n}(\beta_i+1)\ts S_{n,-n}(\beta_i)\ts\eta
\non\end{equation}
with $\gamma$ satisfying \eqref{coga} form a basis of $L$.
\eco

\subsection{Yangian action on the multiplicity space}\label{subsec:ya}

Now we construct an algebra homomorphism $\Y^{\pm}(2)\to\Z(\g_n,\g_{n-1})$
and then use it to define an action of $\Y^{\pm}(2)$ on the multiplicity
space $V(\la)^+_{\mu}$.

For $a,b\in\{-n,n\}$ and a complex parameter $u$
introduce the elements  $Z_{ab}(u)$ of the Mickelsson--Zhelobenko algebra
$\Z(\g_n,\g_{n-1})$ by
\beql{bZab}
Z_{ab}(u)=-\Big(\delta_{ab}(u+\rho_n+\frac12)+F_{ab}\Big)
\prod_{i=-n+1}^{n-1}(u+g_i)
+\sum_{i=-n+1}^{n-1}z_{ai}z_{ib}
\prod_{j=-n+1,\ts j\ne  i}^{n-1}\frac{u+g_j}{g_i-g_j}
\end{equation}
in the $B$ case,
\beql{cZab}
Z_{ab}(u)=\Big(\delta_{ab}(u+\rho_n+\frac12)+F_{ab}\Big)
\prod_{i=-n+1}^{n-1}(u+g_i)
-\sum_{i=-n+1}^{n-1}z_{ai}z_{ib}
\prod_{j=-n+1,\ts j\ne  i}^{n-1}\frac{u+g_j}{g_i-g_j}
\end{equation}
in the $C$ case, and
\begin{multline}\label{dZab}
Z_{ab}(u)=-\Bigg(\Big(\delta_{ab}(u+\rho_n+\frac12)+F_{ab}\Big)
\prod_{i=-n+1}^{n-1}(u+g_i)\\
-\sum_{i=-n+1}^{n-1}z_{ai}z_{ib}\ts (u+g_{-i})
\prod_{j=-n+1,\ts j\ne \pm i}^{n-1}\frac{u+g_j}{g_i-g_j}\Bigg)
\frac{1}{2u+1}
\end{multline}
in the $D$ case,
where $g_i=f_i+1/2$ for all $i$.
In particular, it can be verified that each $Z_{n,-n}(u)$ coincides
with the corresponding interpolation polynomial given in \eqref{bcinterp} 
or \eqref{dinterp}.

Consider now the three cases separately. We shall assume $\mu_{n}=-\infty$
in the notation below.

\bigskip
\noindent
{\bf B type case\/}. 

\bth\label{thm:bhomid}
{\rm (i)} The mapping
\beql{bmap}
s_{ab}(u)\mapsto - u^{-2n}\ts Z_{ab}(u),\qquad a,b\in\{-n,n\}
\end{equation}
defines an algebra homomorphism $\Y^+(2)\to \Z(\g_n,\g_{n-1})$.

{\rm (ii)} The
$\ts\Y^+(2)$-module $V(\lambda)^+_{\mu}$ defined via the homomorphism
\eqref{bmap}
is isomorphic
to the direct sum of two irreducible submodules,
$V(\lambda)^+_{\mu}\simeq U\oplus U'$,
where
\begin{align}
U&=L(0,\beta_1)\ot L(\alpha_2,\beta_2)\ot\cdots\ot 
L(\alpha_{n},\beta_{n})\ot W(1/2),
\non\\
U'&=L(-1,\beta_1)\ot L(\alpha_2,\beta_2)\ot\cdots\ot 
L(\alpha_{n},\beta_{n})\ot W(1/2),
\non
\end{align}
if the $\lambda_i$ are integers 
{\rm(}it is supposed that $U'=\{0\}$ if $\beta_1=0${\rm)}; or
\begin{align}
U&=L(-1/2,\beta_1)\ot L(\alpha_2,\beta_2)\ot\cdots\ot 
L(\alpha_{n},\beta_{n})\ot W(0),
\non\\
U'&=L(-1/2,\beta_1)\ot L(\alpha_2,\beta_2)\ot\cdots\ot 
L(\alpha_{n},\beta_{n})\ot W(1),
\non
\end{align}
if the $\lambda_i$ are half-integers, and the following notation
is used
\begin{align}
\alpha_i&=\min\{\lambda_{i-1},\mu_{i-1}\}-i+1, &\qquad&i=2,\dots,n,
\non
\\
\beta_i&=\max\{\lambda_{i},\mu_{i}\}-i+1, &\qquad&i=1,\dots,n.
\non
\end{align}
\eth

\bigskip
\noindent
{\bf C type case\/}. 

\bth\label{thm:chomid}
{\rm (i)} The mapping
\beql{cmap}
s_{ab}(u)\mapsto (u+1/2)\ts u^{-2n}\ts Z_{ab}(u),\qquad a,b\in\{-n,n\}
\end{equation}
defines an algebra homomorphism $\Y^-(2)\to \Z(\g_n,\g_{n-1})$.

{\rm (ii)} The
$\ts\Y^-(2)$-module $V(\lambda)^+_{\mu}$ defined via the homomorphism
\eqref{cmap} is irreducible and isomorphic
to the tensor product
\beq
L(\alpha_1,\beta_1)\ot\cdots\ot L(\alpha_n,\beta_n),
\non
\end{equation}
where $\al_1=-1/2$ and
\begin{align}
\alpha_i&=\min\{\lambda_{i-1},\mu_{i-1}\}-i+1/2,&\qquad i&=2,\dots,n,
\non\\
\beta_i&=\max\{\lambda_{i},\mu_{i}\}-i+1/2,&\qquad i&=1,\dots,n.
\non
\end{align}
\eth

\bigskip
\noindent
{\bf D type case\/}. 

\bth\label{thm:dhomid}
{\rm (i)} The mapping
\beql{dmap}
s_{ab}(u)\mapsto - 2\ts u^{-2n+2}\ts Z_{ab}(u),\qquad a,b\in\{-n,n\}
\end{equation}
defines an algebra homomorphism $\Y^+(2)\to \Z(\g_n,\g_{n-1})$.

{\rm (ii)} The
$\ts\Y^+(2)$-module $V(\lambda)^+_{\mu}$ defined via the homomorphism
\eqref{dmap} is irreducible and isomorphic
to the tensor product
\begin{equation}
L(\alpha_1,\beta_1)\ot\cdots\ot L(\alpha_{n-1},\beta_{n-1})
\ot W(-\alpha_0),
\non
\end{equation}
where $\alpha_1=\min\{-|\lambda_{1}|,-|\mu_{1}|\}-1/2$,\quad
$\alpha_0=\alpha_1+|\lambda_1+\mu_1|$,
\begin{align}
\alpha_i&=\min\{\lambda_{i},\mu_{i}\}-i+1/2,&\qquad i&=2,\dots,n-1,
\non
\\
\beta_i&=\max\{\lambda_{i+1},\mu_{i+1}\}-i+1/2,&\qquad i&=1,\dots,n-1.
\non
\end{align}
\eth

\Outline
Part (i) of Theorems~\ref{thm:bhomid}--\ref{thm:dhomid} is verified by
using the composition of homomorphisms
\beq
\Y^{\pm}(2)\to \C_n\to \Z(\g_n,\g_{n-1}),
\non\end{equation}
where $\C_n$ is the centralizer $\U(\g_n)^{\g_{n-1}}$.
The first arrow is the homomorphism provided by 
the centralizer construction (see \cite{mo:cc}, \cite{o:ty})
while the second is the natural projection.

By the results of \cite{m:fd}, every irreducible finite-dimensional
representation of the twisted Yangian is a highest weight representation.
It contains a unique, up to a constant factor, vector which is
annihilated by $s_{-n,n}(u)$ and which is an eigenvector
of $s_{nn}(u)$. The corresponding eigenvalue (the highest weight)
uniquely determines the representation.
The vectors in $V(\la)^+_{\mu}$ annihilated by $s_{-n,n}(u)$
can be explicitly constructed by using the lowering operators.
One of such vectors is given by
\beq
\xi_{\mu}=\prod_{i=1}^{n-1}\Big(z_{ni}^{\max\{\lambda_i,\mu_i\}-\mu_i}
z_{i,-n}^{\max\{\lambda_i,\mu_i\}-\lambda_i}\Big)\ts \xi,
\non\end{equation}
where $\xi$ is the highest vector of $V(\la)$.
This is the only vector in the $C,D$ cases, while in the $B$ case
there is another one defined by
\beq
\xi'_{\mu}=z^{}_{n0}\ts \xi_{\mu}.
\non\end{equation}
Calculating the eigenvalues of these vectors we conclude
that they respectively coincide with the eigenvalues of the tensor
product of the highest vectors of the modules $L(\al_i,\beta_i)$;
see \eqref{tensprxi}. \endproof

\bre\label{rem:bra}
{\rm
Theorems~\ref{thm:bhomid}--\ref{thm:dhomid} can be proved
without using the
branching rules for the reductions $\spa_{2n}\downarrow\spa_{2n-2}$
and $\oa_N\downarrow\oa_{N-2}$. Therefore, the reduction multiplicities
can be found by calculating the dimension of the space $V(\lambda)^+_{\mu}$.
For instance, in the symplectic case, Theorem~\ref{thm:chomid} gives
\beq%
c(\mu)=\prod_{i=1}^n(\al_i-\be_i+1)
\non\end{equation}
which, of course, coincides with the value provided by the $C$ type
branching rule; see Section~\ref{subsec:br}.
}
\ere

While keeping $\la$ and $\mu$ fixed we let $\nu$ run over
the values determined by the branching rules; see Section~\ref{subsec:br}.
Using the homomorphisms of Theorems~\ref{thm:bhomid}--\ref{thm:dhomid}
we conclude from \eqref{STT} and \eqref{orSTT} that
the element $S_{n,-n}(u)$
acts in the representation
$V(\la)^+_{\mu}$ precisely as the operator $-Z_{n,-n}(u)$,
$Z_{n,-n}(u)$, or $-2Z_{n,-n}(u)$
in the $B,C$ or $D$ cases, respectively.
Thus, by Corollaries~\ref{cor:bassgen} and \ref{cor:orbassgen},
the following vectors $\xi_{\nu}$ form a basis of the space $V(\la)^+_{\mu}$,
where
\beq
\xi_{\nu}=z_{n0}^{\sigma}\ts\prod_{i=1}^{n}Z_{n,-n}(\gamma_i-1)\cdots 
Z_{n,-n}(\beta_i+1)Z_{n,-n}(\beta_i)\ts\xi_{\mu}
\non\end{equation}
in the $B$ case,
\beql{cvlm}
\xi_{\nu}=\prod_{i=1}^{n}Z_{n,-n}(\gamma_i-1)\cdots 
Z_{n,-n}(\beta_i+1)Z_{n,-n}(\beta_i)\ts\xi_{\mu}
\end{equation}
in the $C$ case, and
\beq
\xi_{\nu}=\prod_{i=1}^{n-1}Z_{n,-n}(\gamma_i-1)\cdots 
Z_{n,-n}(\beta_i+1)Z_{n,-n}(\beta_i)\ts\xi_{\mu}
\non\end{equation}
in the $D$ case. 
Applying the interpolation properties of the polynomials 
$Z_{n,-n}(u)$ we bring the above formulas to the form
given in Lemmas~\ref{lem:bbasms}, 
\ref{lem:cbasms} and \ref{lem:dbasms}, respectively.
Clearly, 
Theorems~\ref{thm:bbasis}, \ref{thm:cbasis} 
and \ref{thm:dbasis} follow.

\subsection{Calculation of the matrix elements}\label{subsec:cm}

Without writing down all explicit formulas
we shall demonstrate how the matrix elements of the generators of
$\g_n$ in the basis $\xi_{\La}$ 
provided by Theorems~\ref{thm:bbasis}, \ref{thm:cbasis} 
and \ref{thm:dbasis} can be calculated.
The interested reader is refered to the papers \cite{m:br, m:wb, m:wbg}
for details.
We choose the following generators
\beq
F_{k-1,-k},\quad F_{k-1,k},\qquad k=1,\dots,n
\non\end{equation}
in the $B$ case,
\beq
F_{k-1,-k},\qquad k=2,\dots,n,\qquad\text{and}\qquad F_{-k,k},\quad F_{k,-k},
\qquad k=1,\dots,n
\non\end{equation}
in the $C$ case, and
\beq
F_{k-1,-k},\quad F_{k-1,k},\qquad k=2,\dots,n, \qquad\text{and}\qquad
F_{21},\quad F_{-2,1}
\non\end{equation}
in the $D$ case.

In the symplectic case the elements 
$F_{kk}$, $F_{k,-k}$, $F_{-k,k}$ commute
with the subalgebra $\g_{k-1}$ in $\U(\g_k)$. Therefore, these
operators preserve the subspace of $\g_{k-1}$-highest vectors
in $V(\lambda)$. So, it suffices to compute the action
of these operators with $k=n$ in the basis $\{\xi_{\nu}\}$
of the space $V(\lambda)^+_{\mu}$; see Lemma~\ref{lem:cbasms}.
For $F_{nn}$ we immediately get
\beq
F_{nn}\ts\xi_{\nu}=\left(2\ts\sum_{i=1}^n\nu_i-
\sum_{i=1}^n\lambda_i-\sum_{i=1}^{n-1}\mu_i\right)\xi_{\nu}.
\non\end{equation}
Further, by \eqref{cvlm}
\beq
Z_{n,-n}(\gamma_i)\ts\xi_{\nu}=\xi_{\nu+\delta_i},\qquad i=1,\dots,n.
\non\end{equation}
However, $Z_{n,-n}(u)$ is a polynomial in $u^2$ of degree $n-1$
with the highest coefficient $F_{n,-n}$. Applying 
the Lagrange interpolation
formula with the interpolation points 
$\gamma_i$, $i=1,\dots,n$
we obtain
\beq
Z_{n,-n}(u)\ts\xi_{\nu}=\sum_{i=1}^n\prod_{a=1,\ts a\ne i}^n
\frac{u^2-\gamma_a^2}{\gamma_i^2-\gamma_a^2}\ \xi_{\nu+\delta_i}.
\non\end{equation}
Taking here the coefficient at $u^{2n-2}$ we get
\beql{Fn-nxi}
F_{n,-n}\ts\xi_{\nu}=\sum_{i=1}^n\prod_{a=1,\ts a\ne i}^n
\frac{1}{\gamma_i^2-\gamma_a^2}\ \xi_{\nu+\delta_i}.
\end{equation}
The action of $F_{-n,n}$ is found in a similar way
with the use of Theorem~\ref{thm:basy}.

In the orthogonal case the action of $F_{nn}$ 
is found in the same way. However,
the elements $F_{n,-n}$ and $F_{-n,n}$
are zero. We shall use second order elements
of the enveloping algebra instead. These are given by
\beq
\Phi_{-a,a}=\frac12\ts\sum_{i=-n+1}^{n-1}F_{-a,i}F_{ia}
\non\end{equation}
with $a\in\{-n,n\}$.
The elements $\Phi_{-a,a}$ commute with the subalgebra $\g_{n-1}$
so that, like in the symplectic case, they preserve the space $V(\lambda)^+_{\mu}$
and their action in the basis $\{\xi_{\nu}\}$ is given by
formulas similar to those for $F_{-a,a}$.

The calculation of the matrix elements of the generators $F_{k-1,-k}$
is similar in all the three cases.
We may assume $k=n$. The operator $F_{n-1,-n}$ preserves
the subspace of $\g_{n-2}$ highest vectors in $V(\lambda)$.
Consider the symplectic case as an example.
Suppose that $\mu'$ is a fixed 
$\g_{n-2}$ highest weight, $\nu'$ is an $(n-1)$-tuple of integers
such that the inequalities \eqref{cineq}
are satisfied with $\lambda$, $\nu$, $\mu$ respectively replaced
by $\mu$, $\nu'$, $\mu'$, and set $\gamma'_i=\nu'_i+\rho_i+1/2$.
It suffices to calculate the action of $F_{n-1,-n}$
on the basis vectors of the form
\beq
\xi_{\nu\mu\nu'}=X_{\mu\nu'}\ts\xi_{\nu\mu},
\non\end{equation}
where $\xi_{\nu\mu}=\xi_{\nu}$ and $X_{\mu\nu'}$ denotes the operator
\beq
X_{\mu\nu'}=\prod_{i=1}^{n-2}
z_{n-1,i}^{\nu'_i-\mu'_i}\ts z_{n-1,-i}^{\nu'_i-\mu_i}\cdot
\prod_{a=m_{n-1}}^{\gamma'_{n-1}-1}
Z_{n-1,-n+1}(a),
\non\end{equation}
where we have used the notation
$m_i=\mu_i+\rho_i+1/2$.
The operator $F_{n-1,-n}$ is permutable with the 
elements $z_{n-1,i}$ and
$Z_{n-1,-n+1}(u)$. Hence, we can write
\beq
F_{n-1,-n}\ts \xi_{\nu\mu\nu'}=X_{\mu\nu'}\ts F_{n-1,-n}\ts \xi_{\nu\mu}.
\non\end{equation}
Now we apply Lemma~\ref{lem:bcdflow}. It remains to calculate
$z_{ni}\ts \xi_{\nu\mu}$ and $X_{\mu\nu'}\ts z_{n-1,-i}$. Using
the relations between the elements of the Mickelsson--Zhelobenko 
algebra $\Z(\g_n,\g_{n-1})$
given in \eqref{bcdrelcom}, we find that
\beq
z_{ni}\ts \xi_{\nu\mu}=\xi_{\nu,\mu-\delta_i}
\non\end{equation}
if $i>0$. Otherwise, if $i=-j$ with positive $j$, write 
\beql{zniximuneg}
z_{n,-j}\ts \xi_{\nu\mu}=z_{n,-j}z_{nj}\ts\xi_{\nu,\mu+\delta_j}=
Z_{n,-n}(m_j)\ts\xi_{\nu,\mu+\delta_j},
\end{equation}
where we have used the
interpolation properties \eqref{zeval} of the polynomials $Z_{n,-n}(u)$.
Finally, we use the expression \eqref{cvlm} of the basis vectors
and Theorem~\ref{thm:basy} to present \eqref{zniximuneg} as a linear combination
of basis vectors. The same argument applies to calculate 
$X_{\mu\nu'}\ts z_{n-1,-i}$. 

The final formulas for the matrix elements of the generators $F_{n-1,-n}$
in all the three cases are given by multiplicative expressions in the entries
of the patterns which exhibit some similarity 
to the formulas of Theorem~\ref{thm:abasis}.

In the orthogonal case we also need to find the action of the generators $F_{n-1,n}$.
Unlike the case of the generators $F_{n-1,-n}$,
the corresponding matrix elements will be given by certain combinations
of multiplicative expressions which do not seem to be possible
to bring to a product form.
There are two alternative ways to calculate
these combinations which we briefly outline below.
First, as in the previous calculation,
we can write
\beq
F_{n-1,n}\ts \xi_{\nu\mu\nu'}=X_{\mu\nu'}\ts F_{n-1,n}\ts \xi_{\nu\mu}.
\non\end{equation}
Applying again Lemma~\ref{lem:bcdflow}, we come to the calculation
of $z_{in}\ts \xi_{\nu\mu}$. This time the interpolation property of
$Z_{-n,-n}(u)$ (see \eqref{bZab} and \eqref{dZab}) allows us to write, e.g., for $i>0$
\beq
z_{in}\ts \xi_{\nu\mu}=z_{in}z_{ni}\ts\xi_{\nu,\mu+\delta_i}=
z_{-n,-i}z_{-i,-n}\ts\xi_{\nu,\mu+\delta_i}=Z_{-n,-n}(m_i)\ts\xi_{\nu,\mu+\delta_i}.
\non\end{equation}
Now, as $Z_{-n,-n}(u)$ is, up to a multiple, the image of $S_{-n,-n}(u)$
under the homomorphism $\Y^{+}(2)\to\Z(\g_n,\g_{n-1})$,
we can express
this operator in terms of the Yangian operators $T_{ab}(u)$
and then apply Theorem~\ref{thm:basy} to calculate its action.

Alternatively, the generator $F_{n-1,n}$ can be written 
modulo the left ideal $\J'$ of $\U'(\g_n)$ as
\beql{Fphi}
F_{n-1,n}
=\Phi_{n-1,-n}(2)\ts\Phi_{-n,n}-
\Phi_{-n,n}\Phi_{n-1,-n}(0),
\end{equation}
where
\beql{bPhinn}
\Phi_{n-1,-n}(u)=\sum_{i=-n+1}^{n-1}z_{n-1,i}\ts z_{i,-n}
\prod_{a=-n+1,\ts a\ne i}^{n-1}\frac{1}{f_i-f_a}\cdot\frac{1}{u+f_i+F_{nn}}
\end{equation}
in the $B$ case, and
\beql{dPhinn}
\Phi_{n-1,-n}(u)=\sum_{i=-n+1}^{n-1}z_{n-1,i}\ts z_{i,-n}
\prod_{a=-n+1,\ts a\ne \pm i}^{n-1}\frac{1}{f_i-f_a}\cdot\frac{1}{u+f_i+F_{nn}}
\end{equation}
in the $D$ case.
The action of $\Phi_{n-1,-n}(u)$
is found exactly as that of $F_{n-1,-n}$ and the matrix elements have
a similar multiplicative form.
Note, however, that the formula \eqref{Fphi}, 
regarded as the equality of operators acting on
$V(\la)^+$, 
is only valid provided the denominators
in \eqref{bPhinn} or \eqref{dPhinn} do not vanish.
Therefore, in order to use \eqref{Fphi}, we first consider $V(\la)$
with `generic' entries of $\la$ and 
calculate the matrix elements of $F_{n-1,n}$
as functions
in the entries of the patterns $\La$.
The final explicit formulas can be written in a singularity-free form
and they are valid
in the general case.

\subsection*{Bibliographical notes}\label{subsec:bibbcd}


The exposition here is based upon the author's papers \cite{m:br, m:wb, m:wbg}.
Slight changes in the notation were made in order to present the results
in a uniform manner for all the three cases.
The branching rules for all classical reductions $\oa_N\downarrow\oa_{N-1}$
and $\spa_{2n}\downarrow\spa_{2n-2}$ are due to Zhelobenko~\cite{z:cg}; see
also Hegerfeldt~\cite{h:bt}, King~\cite{k:wm}, Proctor~\cite{p:yt}, 
Okounkov~\cite{o:mn}, Goodman--Wallach~\cite{gw:ri}.
The lowering operators for the symplectic Lie algebras were first
constructed by Mickelsson~\cite{m:lo}; see also Bincer~\cite{b:ml}.
The explicit relations in the algebra $\Z(\spa_{2n},\spa_{2n-2})$
were calculated by Zhelobenko~\cite{z:za}.

The algebra $\Y(n)$ was first studied in the work of Faddeev
and the St.-Petersburg school in relation with the inverse
scattering method; see for instance
Takhtajan--Faddeev~\cite{tf:qi}, Kulish--Sklyanin~\cite{ks:qs}.
The term ``Yangian" was introduced by Drinfeld in \cite{d:ha}.
In that paper he defined the Yangian $\Y(\agot)$
for each simple finite-dimensional Lie
algebra $\agot$. Finite-dimensional irreducible representations
of $\Y(\agot)$ were classified by Drinfeld~\cite{d:nr} with the use
of a previous work by Tarasov~\cite{t:sq, t:im}. 
Theorem~\ref{thm:basy} goes back to this work of Tarasov;
see also \cite{m:gt}, \cite{nt:ry}.	The criterion
of Theorem~\ref{thm:irrcr} is due to Chari and Pressley~\cite{cp:yr}.
It can also be deduced from the results of \cite{t:sq, t:im};
see \cite{m:fd}. The twisted Yangians were introduced by Olshanski~\cite{o:ty};
see also~\cite{mno:yc}. Their finite-dimensional irreducible
representations were classified in the author's paper~\cite{m:fd}
which, in particular, contains the criteria of Theorems~\ref{thm:symirr} and 
\ref{thm:orthirr}. For more details on the (twisted) Yangians
and their applications 
in the classical representation theory 
see the expository papers~\cite{mno:yc}, \cite{m:ya}
and the recent work of Nazarov~\cite{n:rt, n:ry} where, in particular,
the skew representations of the twisted Yangians were studied.  

In some particular cases, bases in $V(\la)$ were
constructed, e.g., by Wong and Yeh~\cite{wy:md},
Smirnov and Tolstoy~\cite{st:np}.

Weight bases for the fundamental representations of $\oa_{2n+1}$
and $\spa_{2n}$ were independently constructed by Donnelly~\cite{d:ec, d:eco, d:ep}
in a different way. He also demonstrated that the bases of his
coincide with those of Theorems~\ref{thm:bbasis} and \ref{thm:cbasis},
up to a diagonal equivalence.

Harada~\cite{h:sg} employed the results of \cite{m:br}
to construct a new integrable (Gelfand--Tsetlin) system
on the coadjoint orbits of the symplectic groups.
This provides an analog of the Guillemin--Sternberg 
construction~\cite{gs:gc}
for the unitary groups.


\section{Gelfand--Tsetlin bases for representations of $\oa_N$}\label{sec:gtoN}
\setcounter{equation}{0}

In this section we sketch the construction of the bases 
proposed originally by Gelfand and Tsetlin in \cite{gt:fdo}.
It is based upon the fact that the restriction 
$\oa_N\downarrow \oa_{N-1}$ is multiplicity-free.
This makes the construction similar to the $\gl_n$ case.
We shall be applying the general method of Mickelsson algebras
outlined in Section~\ref{subsec:ma}. In particular, the 
corresponding branching rules can be derived from Theorem~\ref{thm:mick};
cf. Section~\ref{subsec:ama}.

It will be convenient to change the notation for
the elements of the orthogonal Lie algebra $\oa_N$
used in Section~\ref{sec:osp}.
We shall now use the standard enumeration
of the rows and columns of $N\times N$-matrices
by the numbers $\{1,\dots,N\}$. Define the involution of
this set of indices by setting $i^{\ts\prime}=N-i+1$.
The Lie algebra $\oa_N$ is spanned by the elements
\beql{newf}
F_{ij}=E_{ij}-E_{j^{\ts\prime}i^{\ts\prime}},\qquad i,j=1,\dots,N.
\end{equation}
We shall keep the notation $\g_n$ for $\oa_N$
with $N=2n+1$ or $N=2n$.

The finite-dimensional irreducible representations of $\g_n$ 
are now parametrized by $n$-tuples
$\lambda=(\lambda_1,\dots,\lambda_n)$ 
where the numbers $\lambda_i$ satisfy the conditions
\beq\label{newlamco}
\lambda_i-\lambda_{i+1}\in \ZZ_+ \qquad\text{for}\quad i=1,\dots,n-1,
\end{equation}
and
\beq\label{newlam1}
\bsp
2\ts \lambda_n&\in\ZZ_+ \qquad\text{for}\quad \g_n=\oa_{2n+1},\\
\lambda_{n-1}+\lambda_n&\in\ZZ_+ \qquad\text{for}\quad \g_n=\oa_{2n}.
\end{split}
\end{equation}
Such an $n$-tuple $\lambda$ is called the {\it highest weight\/}
of the corresponding representation which
we shall denote by $V(\lambda)$.
It contains a unique, up to a constant factor, nonzero vector $\xi$
(the {\it highest vector\/}) such that
\beq
\begin{aligned}
F_{ii}\ts\xi&=\lambda_i\ts\xi\qquad
&\text{for}&\quad i=1,\dots,n,\\
F_{ij}\ts\xi&=0\qquad
&\text{for}&\quad 1\leq i<j\leq N.
\end{aligned}
\non\end{equation}

\subsection{Lowering operators for the reduction 
$\oa_{2n+1}\downarrow\oa_{2n}$}\label{subsec:bdma}

Taking $N=2n+1$ in the definition \eqref{newf}, we shall consider
$\oa_{2n}$ as the subalgebra of $\oa_{2n+1}$ spanned by the elements
\eqref{newf} with $i,j\ne n+1$.
In accordance with the branching rule, the restriction of $V(\la)$
to the subalgebra $\oa_{2n}$ is given by
\beq
V(\lambda)|^{}_{\oa_{2n}}\simeq\underset{\mu}\oplus \ts V'(\mu),
\non\end{equation}
where $V'(\mu)$ is the irreducible finite-dimensional representation of
$\oa_{2n}$ with the highest weight $\mu$ and the sum is taken over
the weights $\mu$ satisfying the inequalities
\beql{bdineq}
\lambda_1\geq\mu_1\geq\lambda_2\geq\mu_2\geq\cdots\geq
\lambda_{n-1}\geq \mu_{n-1}\geq\lambda_{n}\geq|\mu_n|,
\end{equation}
with all the $\mu_i$ being simultaneously
integers or half-integers together with the $\lambda_i$.

The elements $F_{n+1,i}$
span the $\oa_{2n}$-invariant complement
to $\oa_{2n}$ in $\oa_{2n+1}$. Therefore, by the general theory
of Section~\ref{subsec:ma}, the Mickelsson--Zhelobenko algebra $\Z(\oa_{2n+1},\oa_{2n})$
is generated by the elements
\beql{selema}
pF_{n+1,i},\qquad i=1,\dots,n,n^{\ts\prime},\dots,1^{\ts\prime},
\end{equation}
where $p$ is the extremal projector for the Lie algebra  $\oa_{2n}$.
Let $\{\ve_1,\dots,\ve_n\}$ be the basis of $\h^*$ dual to
the basis $\{F_{11},\dots,F_{nn}\}$ of the Cartan subalgebra $\h$ of
$\oa_{2n}$. Set $\ve_{i^{\ts\prime}}=-\ve_i$ for $i=1,\dots,n$.
Denote by $p_{ij}$ the element $p_{\al}$ given by
\eqref{palpha} for the positive root $\al=\ve_i-\ve_j$.
Choosing an appropriate normal ordering on the positive roots,
for any $i=1,\dots,n$ we can write the elements \eqref{selema} in the form
\beql{newselema}
pF_{n+1,i}=	p_{i,i+1}\cdots p_{in}\ts p_{in'}\cdots p_{i1'}\ts F_{n+1,i},
\end{equation}
where the factor $p_{ii^{\ts\prime}}$ is skipped in the product.  
Therefore the right denominator of
this fraction is
\beq
\pi_i=f_{i,i+1}\cdots f_{in}\ts f_{in'}\cdots f_{i1'},
\non\end{equation}
where
\beq
f_{ij}=\begin{cases} F_{ii}-F_{jj}+j-i \qquad&\text{if}\quad j=1,\dots,n\\
F_{ii}-F_{jj}+j-i-2 \qquad&\text{if}\quad j=1',\dots,n'.
\end{cases}
\non\end{equation}
Hence, the elements $s^{\ts\prime}_{ni}=	pF_{n+1,i}\ts \pi_i$ with $i=1,\dots,n$
belong to the Mickelsson algebra $\Sr(\oa_{2n+1},\oa_{2n})$.
One can verify that they are pairwise commuting.

Denote by $V(\lambda)^{+}$ the subspace of $\oa_{2n}$-highest vectors
in $V(\lambda)$.
Given a $\oa_{2n}$-highest weight
$\mu=(\mu_1,\dots,\mu_{n})$ we denote by $V(\lambda)^{+}_{\mu}$
the corresponding weight subspace in $V(\lambda)^{+}$:
\begin{equation}
V(\lambda)^{+}_{\mu}=\{\eta\in V(\lambda)^{+}\ |\ F_{ii}\ts\eta=
\mu_i\ts\eta,\qquad i=1,\dots,n\}.
\non
\end{equation}

By the branching rule, the space $V(\lambda)^{+}_{\mu}$ is 
one-dimensional if the condition \eqref{bdineq} is satisfied.
Otherwise, it is zero.

\bth\label{thm:brabd}
Suppose that the inequalities  \eqref{bdineq} hold.
Then the space $V(\lambda)^{+}_{\mu}$ is spanned
by the vector
\beq
s_{n1}^{\ts\prime\ts\la_1-\mu_1}\cdots s_{nn}^{\ts\prime\ts\la_n-\mu_n}\ts \xi.
\non\end{equation}
\eth

\subsection{Lowering operators for the reduction 
$\oa_{2n}\downarrow\oa_{2n-1}$}\label{subsec:dbma}

Taking $N=2n$ in the definition \eqref{newf}, we shall consider
$\oa_{2n-1}$ as the subalgebra of $\oa_{2n}$ spanned by the elements
\eqref{newf} with $i,j\ne n,n'$ together with
\beq
\frac{1}{\sqrt 2}\ts(F_{ni}-F_{n'i}),
\qquad i=1,\dots,n-1,(n-1)',\dots,1'.
\non\end{equation}

In accordance with the branching rule, the restriction of $V(\la)$
to the subalgebra $\oa_{2n-1}$ is given by
\beq
V(\lambda)|^{}_{\oa_{2n-1}}\simeq\underset{\mu}\oplus \ts V'(\mu),
\non\end{equation}
where $V'(\mu)$ is the irreducible finite-dimensional representation of
$\oa_{2n-1}$ with the highest weight $\mu$ and the sum is taken over
the weights $\mu$ satisfying the inequalities
\beql{ddbineq}
\lambda_1\geq\mu_1\geq\lambda_2\geq\mu_2\geq\cdots\geq\lambda_{n-1} \geq
\mu_{n-1}\geq|\lambda_{n}|,
\end{equation}
with all the $\mu_i$ being simultaneously
integers or half-integers together with the $\lambda_i$.

The elements
\beql{bdcompnew}
F_{nn},\qquad
F^{\ts\prime}_{ni}=\frac{1}{\sqrt 2}\ts(F_{ni}+F_{n'i}),\qquad
i=1,\dots,n-1,(n-1)',\dots,1'
\end{equation}
span the $\oa_{2n-1}$-invariant complement
to $\oa_{2n-1}$ in $\oa_{2n}$. Therefore, by the general theory
of Section~\ref{subsec:ma}, the Mickelsson--Zhelobenko algebra $\Z(\oa_{2n},\oa_{2n-1})$
is generated by the elements
\beql{dselema}
pF_{nn},\qquad pF^{\ts\prime}_{ni},\qquad i=1,\dots,n-1,(n-1)',\dots,1',
\end{equation}
where $p$ is the extremal projector for the Lie algebra  $\oa_{2n-1}$.
Let $\{\ve_1,\dots,\ve_{n-1}\}$ be the basis of $\h^*$ dual to
the basis $\{F_{11},\dots,F_{n-1,n-1}\}$ of the Cartan subalgebra $\h$ of
$\oa_{2n-1}$. Set $\ve_{i^{\ts\prime}}=-\ve_i$ for $i=1,\dots,n-1$.
Denote by $p_{ij}$ and $p_i$ the elements $p_{\al}$ given by
\eqref{palpha} for the positive roots $\al=\ve_i-\ve_j$
and $\al=\ve_i$, respectively.
Choosing an appropriate normal ordering on the positive roots,
for any $i=1,\dots,n-1$
we can write the elements \eqref{dselema} in the form
\beql{dnewselema}
pF^{\ts\prime}_{ni}=	p_{i,i+1}\cdots p_{i,n-1}\ts p_i\ts 
p_{i,(n-1)'}\cdots p_{i1'}\ts F^{\ts\prime}_{ni},
\end{equation}
where the factor $p_{ii^{\ts\prime}}$
is skipped in the product.
Therefore the right denominator of
this fraction is
\beq
\pi_i=f^{}_{i,i+1}\cdots f^{}_{i,n-1}\ts
f^{}_{i}\ts f'_{i}\ts f^{}_{i,(n-1)'}\cdots f^{}_{i1'},
\non\end{equation}
where
\beq
f^{}_{ij}=\begin{cases} F_{ii}-F_{jj}+j-i \qquad&\text{if}\quad j=1,\dots,n-1\\
F_{ii}-F_{jj}+j-i-2 \qquad&\text{if}\quad j=1',\dots,(n-1)'
\end{cases}
\non\end{equation}
and $f^{}_{i}=f'_{i}-1=2(F_{ii}+n-i)$.
Hence, the elements $s_{ni}=pF^{\ts\prime}_{ni}\ts \pi_i$ with $i=1,\dots,n-1$
belong to the Mickelsson algebra $\Sr(\oa_{2n},\oa_{2n-1})$.
One can verify that they are pairwise commuting.

Denote by $V(\lambda)^{+}$ the subspace of $\oa_{2n-1}$-highest vectors
in $V(\lambda)$.
Given a $\oa_{2n-1}$-highest weight
$\mu=(\mu_1,\dots,\mu_{n-1})$ we denote by $V(\lambda)^{+}_{\mu}$
the corresponding weight subspace in $V(\lambda)^{+}$:
\begin{equation}
V(\lambda)^{+}_{\mu}=\{\eta\in V(\lambda)^{+}\ |\ F_{ii}\ts\eta=
\mu_i\ts\eta,\qquad i=1,\dots,n-1\}.
\non
\end{equation}

By the branching rule, the space $V(\lambda)^{+}_{\mu}$ is 
one-dimensional if the condition \eqref{ddbineq} is satisfied.
Otherwise, it is zero.

\bth\label{thm:drabd}
Suppose that the inequalities  \eqref{ddbineq} hold.
Then the space $V(\lambda)^{+}_{\mu}$ is spanned
by the vector
\beq
s_{n1}^{\la_1-\mu_1}\cdots s_{n,n-1}^{\la_{n-1}-\mu_{n-1}}\ts \xi.
\non\end{equation}
\eth

Note that the generator $pF_{nn}$ of the algebra $\Z(\oa_{2n},\oa_{2n-1})$
does not occur in the formula for the basis vector as it has the zero weight
with respect to $\h$.

\subsection{Basis vectors}\label{subsec:bdbv}

The representation $V(\lambda)$ of the Lie algebra $\g_n=\oa_{2n+1}$ 
or $\oa_{2n}$
is equipped with a contravariant inner product which is 
uniquely determined by the conditions
\beq
\lan\xi,\xi\ran=1\qquad\text{and}\qquad  \lan F_{ij}\ts u,v\ran=\lan u,F_{ji}\ts v\ran
\non\end{equation}
for all $u,v\in V(\la)$ and any indices $i,j$.

Combining Theorems~\ref{thm:brabd} and \ref{thm:drabd} we can construct
another basis for each representation $V(\lambda)$ of $\g_n$; cf. Section~\ref{subsec:br}.

\bigskip
\noindent
{\bf B type case\/}. We need to modify the definition
of the $B$ type 
pattern $\Lambda$ introduced in Section~\ref{subsec:br}. 
Here 
$\Lambda$ is an array of the form
\begin{align}
\quad&\lambda^{}_{n1}\qquad\lambda^{}_{n2}
\qquad\qquad\cdots\qquad\qquad\lambda^{}_{nn}\non\\
&\qquad\lambda'_{n1}\qquad  \lambda'_{n2}
\qquad\qquad\cdots\qquad\qquad\lambda'_{nn}\non\\
\ &\ \ \qquad\qquad\lambda^{}_{n-1,1}\qquad\cdots
\qquad\lambda^{}_{n-1,n-1}\non\\
&\quad\qquad\qquad\qquad\lambda'_{n-1,1}
\qquad\cdots\qquad\lambda'_{n-1,n-1}\non\\
&\qquad\qquad\qquad\qquad\qquad\cdots\qquad\cdots\non\\
\quad&\ \ \qquad\qquad\qquad\qquad\qquad\qquad\qquad\lambda^{}_{11}\non\\
&\ \ \qquad\qquad\qquad\qquad\qquad\qquad\qquad\qquad\lambda'_{11}\non
\end{align}
such that $\lambda=(\lambda^{}_{n1},\dots, \lambda^{}_{nn})$,
the remaining
entries are all
integers or half-integers together with the
$\lambda_i$, and the following inequalities hold
\beq
\lambda^{}_{k1}\geq\lambda'_{k1}\geq\lambda^{}_{k2}\geq\lambda'_{k2}\geq
\cdots\geq
\lambda'_{k,k-1}\geq\lambda^{}_{kk}\geq|\lambda'_{kk}|
\non
\end{equation}
for $k=1,\dots,n$, and
\beq
\lambda'_{k1}\geq\lambda^{}_{k-1,1}\geq\lambda'_{k2}\geq
\lambda^{}_{k-1,2}\geq \cdots\geq
\lambda'_{k,k-1}\geq\lambda^{}_{k-1,k-1}\geq|\lambda'_{kk}|
\non
\end{equation}
for $k=2,\dots,n$.

\bth\label{thm:bbasiso}
The vectors
\beq
\eta^{}_{\Lambda}=s_{11}^{\ts\prime\ts\la^{}_{11}-\la'_{11}}\ts
\prod_{k=2,\dots,n}^{\rightarrow}
\Big(s_{k1}^{\ts\prime\ts \la^{}_{k1}-\la'_{k1}}\ts
\cdots s_{kk}^{\ts\prime\ts \la^{}_{kk}-\la'_{kk}} 
s_{k1}^{\la'_{k1}-\la^{}_{k-1,1}}\cdots 
s_{k,k-1}^{\la'_{k,k-1}-\la^{}_{k-1,k-1}}\Big)\ts\xi
\non\end{equation}
parametrized by the 
patterns $\Lambda$ form an orthogonal basis of 
the representation $V(\lambda)$. 
\eth

\bigskip
\noindent
{\bf D type case\/}. 	Here we define
the $D$ type 
patterns $\Lambda$ as arrays of the form
\begin{align}
&\qquad\lambda^{}_{n1}\qquad  \lambda^{}_{n2}
\qquad\qquad\cdots\qquad\qquad\lambda^{}_{nn}\non\\
\ &\ \ \qquad\qquad\lambda'_{n-1,1}\qquad\cdots
\qquad\lambda'_{n-1,n-1}\non\\
&\quad\qquad\qquad\qquad\lambda^{}_{n-1,1}
\qquad\cdots\qquad\lambda^{}_{n-1,n-1}\non\\
&\qquad\qquad\qquad\qquad\qquad\cdots\qquad\cdots\non\\
\quad&\qquad\qquad\qquad\qquad\qquad\qquad\qquad\lambda'_{11}\non\\
&\qquad\qquad\qquad\qquad\qquad\qquad\qquad\qquad\lambda^{}_{11}\non
\end{align}
such that $\lambda=(\lambda^{}_{n1},\dots, \lambda^{}_{nn})$,
the remaining
entries are all
integers or half-integers together with the
$\lambda_i$, and the following inequalities hold
\beq
\lambda^{}_{k1}\geq\lambda'_{k-1,1}\geq\lambda^{}_{k2}\geq
\lambda'_{k-1,2}\geq \cdots\geq
\lambda^{}_{k,k-1}\geq\lambda'_{k-1,k-1}\geq|\lambda^{}_{kk}|
\non
\end{equation}
for $k=2,\dots,n$, and 
\beq
\lambda'_{k1}\geq\lambda^{}_{k1}\geq\lambda'_{k2}\geq\lambda^{}_{k2}\geq
\cdots\geq
\lambda^{}_{k,k-1}\geq\lambda'_{kk}\geq|\lambda^{}_{kk}|
\non
\end{equation}
for $k=1,\dots,n-1$.

\bth\label{thm:dbasiso}
The vectors
\beq
\eta^{}_{\Lambda}=
\prod_{k=1,\dots,n-1}^{\rightarrow}
\Big(s_{k+1,1}^{\la^{}_{k+1,1}-\la'_{k1}}\ts
\cdots s_{k+1,k}^{\la^{}_{k+1,k}-\la'_{kk}} 
s_{k1}^{\ts\prime\ts\la'_{k1}-\la^{}_{k1}}\cdots 
s_{kk}^{\ts\prime\ts\la'_{kk}-\la^{}_{kk}}\Big)\ts\xi
\non\end{equation}
parametrized by the 
patterns $\Lambda$ form an orthogonal basis of 
the representation $V(\lambda)$. 
\eth

The norms of the basis vectors
$\eta^{}_{\Lambda}$ can be found in an explicit form.
The formulas for the matrix
elements of the generators of the Lie algebra
$\oa_{N}$ in the original paper by Gelfand and Tsetlin~\cite{gt:fdo} 
are given in the orthonormal basis
\beq
\zeta^{}_{\Lambda}=\eta^{}_{\Lambda}/ \| \eta^{}_{\Lambda}\|,\qquad
\| \eta^{}_{\Lambda}\|^2=\lan\eta^{}_{\Lambda},	\eta^{}_{\Lambda}\ran.
\non\end{equation}

\subsection*{Bibliographical notes}\label{subsec:bibbd}


The exposition of this section follows Zhelobenko~\cite{z:gz}.
The branching rules were previously derived by him in \cite{z:cg}.
The lowering operators for the reduction
$\oa_N\downarrow \oa_{N-1}$ were constructed by Pang and Hecht~\cite{ph:lr} 
and Wong~\cite{w:ro}; see also Mickelsson~\cite{m:lor}.
They are presented in a form similar to \eqref{zin}
and \eqref{zni} although more complicated. A derivation of the matrix
element formulas of Gelfand and Tsetlin~\cite{gt:fdo} was also 
given in \cite{ph:lr} and \cite{w:ro} which basically
follows the approach outlined in Section~\ref{subsec:loexp}.
The defining relations for the algebra $\Z(\oa_N,\oa_{N-1})$
were given in an explicit form by Zhelobenko~\cite{z:za}.
Gould's approach based upon the characteristic identities of
Bracken and Green~\cite{bg:vo, gr:ci} for the orthogonal Lie algebras
is also applicable; see Gould~\cite{g:cir, g:ia, g:wc}. It produces an independent
derivation of the matrix element formulas. Although the quantum minor approach
has not been developed so far for the Gelfand--Tsetlin basis
for the orthogonal Lie algebras,
it seems to be plausible that the corresponding analogs of
the results outlined in Section~\ref{subsec:qm}
can be obtained.

Analogs of the Gelfand--Tsetlin bases \cite{gt:fdo} for 
representations of a nonstandard deformation $\U'_q(\oa_N)$ of $\U(\oa_N)$
were given by Gavrilik and Klimyk~\cite{gk:qd},
Gavrilik and Iorgov~\cite{gi:qd} and Iorgov and Klimyk~\cite{ik:nd}.

The Gelfand--Tsetlin modules over the orthogonal Lie algebras
were studied by Mazorchuk~\cite{m:gz} with the use of the
matrix element formulas from \cite{gt:fdo}.



\end{document}